\documentclass[a4paper]{siamart1116}
\usepackage{algorithm,algorithmic}
\usepackage{amsmath}
\usepackage{amssymb}
\graphicspath{{images/}}

\title{Parallelization and scalability analysis of inverse factorization using the Chunks and Tasks programming model\thanks{\today
\funding{This work was supported by the Swedish research council under grant 621-2012-3861 and the Swedish national strategic e-science research program (eSSENCE).}}}

\author{Anton G. Artemov%
\thanks{Division of Scientific Computing, Department of Information Technology, Uppsala University, Box 337, SE-751 05 Uppsala, Sweden (\email{anton.artemov@it.uu.se},  \email{elias.rudberg@it.uu.se}, \email{emanuel.rubensson@it.uu.se}).}
\and
Elias Rudberg%
\footnotemark[2]
\and
Emanuel H. Rubensson%
\footnotemark[2]
}

\begin{document}
\maketitle

\begin{abstract}
We present three methods for distributed memory parallel inverse factorization of block-sparse Hermitian positive definite matrices. The three methods are a recursive variant of the AINV inverse Cholesky algorithm, iterative refinement, and localized inverse factorization, respectively. All three methods are implemented using the Chunks and Tasks programming model, building on the distributed sparse quad-tree matrix representation and parallel matrix-matrix multiplication in the publicly available Chunks and Tasks Matrix Library (CHTML). Although the algorithms are generally applicable, this work was mainly motivated by the need for efficient and scalable inverse factorization of the basis set overlap matrix in large scale electronic structure calculations. We perform various computational tests on overlap matrices for quasi-linear Glutamic Acid-Alanine molecules and three-dimensional water clusters discretized using the standard Gaussian basis set STO-3G with up to more than 10 million basis functions. We show that for such matrices the computational cost increases only linearly with system size for all the three methods. We show both theoretically and in numerical experiments that the methods based on iterative refinement and localized inverse factorization outperform previous parallel implementations in weak scaling tests where the system size is increased in direct proportion to the number of processes. We show also that compared to the method based on pure iterative refinement the localized inverse factorization requires much less communication.
\end{abstract}

\section{Introduction}\label{sec:introduction_new}
We focus in this work on the efficient and scalable computation of
inverse factors of symmetric positive definite matrices.  Although the
methods considered are generally applicable, we are in particular
interested in the factorization of basis set overlap matrices coming
from discretizations using localized basis sets in large scale
electronic structure calculations based on, for instance,
Hartree--Fock~\cite{levine2009quantum} or Kohn--Sham density
functional theory~\cite{kohn1965self}. In this context, the inverse
factor is used in a congruence transformation of the
Roothaan--Hall/Kohn--Sham eigenvalue equations:
\begin{equation}
  F x = \lambda S x \qquad \longrightarrow \qquad (Z^*FZ)y = \lambda y
\end{equation}
where $F$ is the Fock or Kohn--Sham matrix, $S$ is the overlap matrix,
and $S^{-1} = ZZ^*$ is an inverse factorization thereof.

Electronic structure theory is of fundamental importance in a number
of scientific disciplines, including biology, chemistry, and materials
science. Unfortunately, the solution of the quantum mechanical
eigenvalue problem outlined above using standard methods is extremely
expensive and usually limited to rather small molecular systems.  To
be able to reach larger molecular systems and materials there are
essentially two possibilities: reduce the computational workload or
increase the computational power. The workload of standard electronic
structure methods increases cubically with system size. This means
that an increase of the computational power alone has a very limited
effect on the system sizes that are within reach. Therefore, so-called
linear scaling methods have been developed for which the workload
increases only linearly with system size. A remarkable progress of
such methods has occurred in the last 15-20 years,
see~\cite{Bowler_2012} for an excellent overview of this
development. With a linear scaling method at hand it makes sense to
attempt to reach larger systems by increasing the computational power.
Today this is possible only with a parallel implementation. Ideally,
the runtime should stay constant as the workload is increased in
direct proportion to the number of processors; this is sometimes referred to
as linear weak scaling efficiency. With both linear scaling of the
total workload with system size and linear weak scaling efficiency,
the system sizes within reach should increase linearly with the
computational power employed.  Therefore, our primary criteria for
assessment of the developed methods are scaling with system size and
weak scaling efficiency. With the current trend that the cost of
floating point operations relative to the cost of movement of data is
becoming smaller we shall also consider the communication costs to be
important.

Two approaches have dominated previous work on inverse factorization of
the overlap matrix. The first one is computation of the inverse
Cholesky factor using the AINV algorithm~\cite{benzi1996sparse}. The
inverse Cholesky factor has frequently been used in electronic
structure calculations~\cite{Millam1997, Challacombe1999, Xiang2007}. 
There are several variants of the AINV algorithm including a
stabilized one~\cite{benzi2000robust}, a blocked
one~\cite{benzi2001stabilized} and a recursive
one~\cite{rubensson2007hierarchic}, where the latter is suitable for
hierarchical matrix representations.
The AINV algorithms are very fast but, as will be shown in the present
work, the parallelism in those algorithms is limited. 
The second approach is computation of the inverse square root using an
iterative refinement method~\cite{niklasson2004iterative} closely
related to the Newton--Schulz iteration for the sign matrix
function~\cite{higham2008functions,localized_inverse_factorization}. This
method is based on sparse matrix-matrix multiplication and has also
been frequently used in electronic structure calculations, using a
scaled identity matrix as starting guess~\cite{jansik2007linear,
  VandeVondele_2012}.  

A third alternative is recursive inverse factorization which is a
divide-and-conquer method that combines iterative refinement with a
recursive decomposition of the matrix~\cite{rubensson2008recursive}.
The matrix is partitioned using a binary principal submatrix
decomposition, corresponding to a cut of the index set into two
pieces.  Inverse factors of the two principal submatrices are computed
and iterative refinement is then used in a correction step to get an
inverse factor of the whole matrix. This procedure is applied
recursively giving the recursive inverse factorization
algorithm. Recently a localized variant of this method was proposed
which uses only localized computations for the correction step with
cost proportional to the cut
size~\cite{localized_inverse_factorization}.  We refer to this
variant as \emph{localized inverse factorization}.  A key feature of
both the regular and localized recursive inverse factorization is that
the inverse factorizations at lower levels are completely disconnected
and therefore embarrassingly parallel.

We consider in this work one variant in each of the three classes of
methods outlined above: 1) the recursive variant of the AINV
algorithm, hereinafter referred to as \emph{the recursive inverse
  Cholesky algorithm} (RINCH), 2) iterative refinement using a scaled identity
as starting guess (IRSI), and 3) localized inverse factorization (LIF).
All three alternatives considered here make use of sparse
matrix-matrix multiplication as a computational kernel.

Sparse matrix-matrix multiplication for general matrices with a priori
unknown sparsity patterns is difficult to parallelize, especially with
good weak scaling performance and taking advantage of data
locality~\cite{communication_optimal_2}. A common approach is to
randomly permute rows and columns of the input matrices to destroy any
structure. Then, some algorithm for dense matrices, such as Cannon's
algorithm or SUMMA, is applied, but with the dense submatrix products
replaced by sparse submatrix products~\cite{BulucGilbert2012,
  Borstnik2014}. In this way load balance is achieved. Linear weak
scaling efficiency is possible if, instead of the two-dimensional
process grid used in Cannon's algorithm and SUMMA, a three-dimensional
process grid is used as
in~\cite{communication_optimal_2,azad20163dmatmul, dawson_ntpoly_2018}. 
However, because of the random permutation of matrix rows and columns,
the possibility to exploit the nonzero structure to avoid movement of
data or communication is lost.

In this article we build on our previous work on a sparse matrix
library based on sparse quaternary tree (quad-tree) matrix
representation implemented using the Chunks and Tasks programming
model~\cite{CHT-PARCO-2014, rubensson2016locality}. The library is
called the Chunks
and Tasks Matrix Library (CHTML) and is publicly available at
\url{chunks-and-tasks.org} under the modified BSD license. This
implementation exploits data locality to achieve efficient
representation of sparse matrices. For sparse matrix-matrix
multiplication, this leads to much less movement of data compared to
approaches that do not take advantage of data locality. In particular,
for matrices with data locality, such as overlap matrices or other
matrices coming from discretizations of physical systems using
localized basis sets, the average amount of data received per process
is constant in the weak scaling scenario outlined above where the
number of processors is scaled up in direct proportion to the system
size. The total execution time is instead usually dominated by the
execution of tasks along the critical path.  For the iterative
refinement this scaling behavior should be directly inherited. 
For the localized inverse factorization the localization in the
computation is even stronger with the bulk of the work lying in
subproblems that are completely disconnected. In theory very little
communication should therefore be needed in a parallel implementation
of the localized inverse factorization algorithm. In practice, using
a traditional message passing approach, it would be very difficult
to achieve a parallel implementation of the localized inverse
factorization algorithm with a balanced distribution of both work and
data for sparse matrices with a priori unknown sparsity patterns.
We show in this article, using instead a Chunks and
Tasks implementation, that reduced communication costs compared to
other scalable algorithms can be realized in practice.

We have implemented the recursive inverse Cholesky algorithm, the
iterative refinement method, and localized inverse factorization in
CHTML. In Section~\ref{sec:algorithms} we describe the three
algorithms. In Section~\ref{sec:cpl_estimation} we analyse the scalability
of the algorithms by considering their critical paths. We are in
particular interested in how the critical paths depend on system size,
since this puts a limit on the weak scaling performance that can be
achieved. In Section~\ref{sec:cht_implementation} we give some background to
the choice of programming model and dicuss some restrictions in the
model and their implications for expressiveness and performance. In
Section~\ref{sec:experiments} we present numerical experiments for overlap
matrices coming from discretizations using standard Gaussian basis
sets, for quasi-linear Glutamic Acid-Alanine molecules as well as
three-dimensional water clusters. It is shown that all three methods
achieve linear scaling with system size for the considered
systems. Weak scaling and actual critical path lengths as reported by
the Chunks and Tasks library are compared to the theoretical results.
We end the article with a discussion and some concluding remarks in
Section \ref{sec:discussion}.
 
\section{Algorithms} \label{sec:algorithms}
In this section we describe the inverse factorization algorithms
considered in this work.  All algorithms are written such as to
operate on matrices represented using quad-trees. 
 
\subsection{The recursive inverse Cholesky algorithm} 
The recursive inverse Cholesky algorithm (RINCH)
\cite{rubensson2007hierarchic} is a variant of the AINV algorithm that
is adapted to hierarchical matrix structures and uses recursive calls
to go through the hierarchy. We use here the left-looking variant of
RINCH given by Algorithm~\ref{RINCH_presudocode}. This algorithm is
equivalent to Algorithm~7 in \cite{rubensson2007hierarchic} with $n=2$
blocks in every dimension. When reaching the leaf level of the
hierarchy, instead of the recursive call, one may for instance use the
regular non-recursive AINV or routines for Cholesky factorization and
inversion of triangular matrices in the Basic Linear Algebra
Subprograms (BLAS) to compute the inverse factor of the leaf matrix.
The most common operation in the algorithm is matrix-matrix
multiplication, which, besides the recursive RINCH calls, is also the
most time-consuming operation.

\begin{algorithm}[t]
\begin{algorithmic}[1]
\REQUIRE Hermitian positive definite matrix $S$
\ENSURE $Z$
\STATE{\textbf{Input} $S$}
\IF{lowest level}
	\STATE{ Factorize $S^{-1} = ZZ^*$ and \textbf{return} $Z$ }
\ELSE
	\STATE{$Z_{0,0} = $ {\sc rinch}$(S_{0,0})$}
	\STATE{$R = Z^*_{0,0} S_{0,1}$} \label{algline:rinch_line6}
	\STATE{$Q = - R^* R + S_{1,1}$} \label{algline:rinch_line7}
	\STATE{$Z_{1,1} = $ {\sc rinch}$(Q)$}
	\STATE{$Z_{0,1} = -Z_{0,0}R Z_{1,1}$} \label{algline:rinch_line9}
	\STATE{$Z_{1,0} = 0$}
	\STATE{$Z =  \begin{pmatrix}
	Z_{0,0} & Z_{0,1} \\ Z_{1,0} & Z_{1,1}
	\end{pmatrix} $}
	\STATE{\textbf{return} $Z$}
\ENDIF

\end{algorithmic}
\caption{Recursive inverse Cholesky, {\sc rinch}$(S)$}
\label{RINCH_presudocode}
\end{algorithm}

\subsection{Iterative refinement}

The refinement procedure proposed in \cite{niklasson2004iterative} can
be used to improve the accuracy of an approximate inverse factor. The
procedure works as follows: given $Z_0$ as an approximation of the
inverse factor of $S$ such that $ \| I - Z^*_0 S Z_0\|_2 < 1,$ a
sequence of matrices $Z_i$ is constructed such that the factorization
error $\delta_i = I - Z^*_i S Z_i \rightarrow 0$ as $i \rightarrow
\infty.$ 
The regular iterative refinement is given by 
\begin{equation}
 b_0 = 1,\,\,  b_k = \frac{2k-1}{2k}b_{k-1},k=1\ldots m,\,\,Z_{i+1} = Z_i \sum \limits_{k=0}^{m}b_k \delta_{i}^{k},\,\, \delta_i = I - Z^{*}_{i}SZ_i. \label{refinement_nonlocal}
\end{equation} 
If it is known that only a small part of the matrix needs to be
updated, then one can use a localized method in order to reduce costs:

\begin{equation}
\begin{cases}
Z_{i+1} = Z_i \sum \limits_{k=0}^{m}b_k \delta_{i}^{k}, \\
\delta_{i+1} = \delta_{i} - Z^{*}_{i+1} S ( Z_{i+1} - Z_{i}) - (Z_{i+1} - Z_{i})^{*} S Z_{i}.
\end{cases} \label{refinement_local}
\end{equation}

One possible choice for the starting guess is a scaled identity matrix
\begin{equation}
Z_0 = cI
\end{equation}
which with an appropriate scaling $c$ gives the inverse square root of
the matrix as a result \cite{higham1997stable, jansik2007linear}.  The
optimal scaling is given by
\begin{equation}
c = \sqrt{\frac{2}{\lambda_{\min} + \lambda_{\max}}}
\end{equation}
where $\lambda_{\min}$ and $\lambda_{\max}$ are the extremal
eigenvalues of $S$ \cite{rubensson2008recursive}.

In our case the exact eigenvalues are not available and we will use an
approximation for the scaling factor suggested in
\cite{jansik2007linear}:

\begin{equation}\label{eq:scaling_gers}
c = \sqrt{\frac{2}{\beta}},
\end{equation} 
where $\beta$ is an approximation of the spectral radius of $S$
computed with the Gershgorin circle theorem \cite{gershgorin1931uber}.

In order to minimize the number of parameters, we employ the automatic
stopping criterion proposed in \cite{kruchinina2016parameterless} and used in \cite{localized_inverse_factorization}. The iterative refinement is
stopped as soon as $\|\delta_{i+1}\|_F > \|\delta_{i}\|_F^{m+1}$.
Using this stopping criterion, the iterative refinement process
automatically terminates when rounding or other numerical errors start
to dominate over the factorization error. 

\subsection{Recursive and localized inverse factorization}
Recursive inverse factorization first described in
\cite{rubensson2008recursive} 
is based on the following result. Let the input matrix $S$ be
partitioned as $S = \begin{bmatrix} A & B \\ B^* & C \end{bmatrix}$
and let $Z_0 = \begin{bmatrix} Z_A & 0 \\ 0 & Z_C \end{bmatrix}$,
where $Z_A$ and $Z_C$ are inverse factors of $A$ and $C$,
respectively. Then, $\|I-Z_0^*SZ_0\|_2 < 1$, which implies convergence
of iterative refinement with $Z_0$ as starting
guess~\cite{rubensson2008recursive}.  This result was recently
strengthened and it was shown that $\|I-Z_0^*SZ_0\|_2 \leq 1 -
\frac{\lambda_{\textrm{min}}(S)}{\lambda_{\textrm{max}}(S)}$~\cite{localized_inverse_factorization}.
Those convergence results immediately suggest a \emph{recursive
  inverse factorization} algorithm where iterative refinement is
combined with a recursive binary principal submatrix decomposition of
the matrix $S$.

Recently we proposed a localized variant of the recursive inverse
factorization algorithm that makes use of the localized iterative
refinement in \eqref{refinement_local} and a local construction of the
starting guess~\cite{localized_inverse_factorization}. An advantage
of this \emph{localized inverse factorization}, given by
Algorithm~\ref{alg:localized_inverse_factorization}, is that under
certain assumptions on $S$ and the recursive principal submatrix
decomposition, the cost of combining $Z_A$ and $Z_C$ to get an inverse
factor of $S$ becomes negligible compared to the overall cost for
sufficiently large systems~\cite{localized_inverse_factorization}.
Another important feature of this algorithm is that $Z_A$ and $Z_C$
can be computed in parallel without any communication between the two
computations.

When the lowest level is reached in the algorithm, leaf-level routines
are called to compute the inverse factor. The leaf level inverse
factorization may for example use BLAS routines or some standard AINV
algorithm. Alternatively, one may switch to for example the RINCH
algorithm before reaching the leaf level of the quad-tree.

\begin{algorithm}[t]
\begin{algorithmic}[1]
\REQUIRE Hermitian positive definite matrix $S$
\ENSURE $Z$
\IF{lowest level}
	\STATE Factorize $S^{-1} = ZZ^*$ and \textbf{return} $Z$
\ENDIF
\STATE Matrix partition $S = \begin{bmatrix} A & B \\ B^* & C \end{bmatrix}$
\STATE $Z_A = $ {\sc lif}$(A)$, $Z_C = $ {\sc lif}$(C)$
\STATE $Z_0 = \begin{bmatrix} Z_A & 0 \\ 0 & Z_C \end{bmatrix}$ 
\STATE $\delta_0 = -\begin{bmatrix} 0 & Z_A^*BZ_C \\ Z_C^*B^*Z_A & 0 \end{bmatrix}$ \label{algline:delta_zero_construction}
\STATE $i = 0$
\REPEAT
\STATE $M_i = Z_i (\sum \limits_{k=1}^{m}b_k \delta_{i}^{k})$ \label{algline:iterrefine_start}
\STATE $Z_{i+1} = Z_i + M_i$ \label{algline:Z_update}
\STATE $\delta_{i+1} = \delta_{i} - Z^{*}_{i+1} (SM_i) - (M_i^{*} S) Z_{i}$ \label{algline:delta_update}
%\STATE $Z_{i+1} = Z_i \sum \limits_{k=0}^{m}b_k \delta_{i}^{k}$
%\STATE $\delta_{i+1} = \delta_{i} - Z^{*}_{i+1} S ( Z_{i+1} - Z_{i}) - (Z_{i+1} - Z_{i})^{*} S Z_{i}$
\STATE $i = i+1$ \label{algline:iterrefine_end}
\UNTIL{$\| \delta_{i+1} \|_F > \| \delta_i \|_F^{m+1}$}
\RETURN $Z_{i+1}$
\end{algorithmic}
\caption{Localized inverse factorization, {\sc lif($S$)}}
\label{alg:localized_inverse_factorization}
\end{algorithm}

\section{Estimation of critical path lengths}\label{sec:cpl_estimation}

For any given system, the computational time will be bounded from below by a value determined by the critical path length. Even if one assumes that the runtime is able to exploit all parallelism in the algorithm and that an infinite number of processes are available, the computational time cannot fall below this value. In this section we will derive critical path lengths as a function of system size for each of the three algorithms considered in this work. In this way we get the best possible weak scaling performance for each of the three algorithms.

First of all, let us define a task as a large enough computation,
which may consist of child tasks or arithmetical
operations, or both. For example, matrix addition is a task and,
taking into account matrix representation as a quad-tree, it actually
consists of 4 child tasks, which are also matrix additions. Each of
those also consists of 4 tasks unless the leaf level of the hierarchy
is reached where arithmetical routines are performed. We
assume that a task can have multiple inputs and a single output.
 
Using the definition of a task, the whole computation can be viewed as
the traversal of a directed graph, where each vertex is a task and
the edges represent data dependencies between tasks.
Then, the critical path can be defined as the longest
path in the directed graph of the computation, which connects start
and end vertices. The critical path is the longest series of
operations which has to be executed serially due to data dependencies
and its length is an important feature of any algorithm from a
parallel programmer's point of view.

In our derivations, we consider matrix-matrix multiplication,
matrix-scalar multiplication, matrix addition, separate matrix transpose and
factorization routines as tasks. Small
computations like multiplication or addition of real numbers are not
considered as tasks. If a matrix transpose is involved in a multiplication, we assume 
that it does not require a separate task.

Let us derive the critical path estimation for the algorithms
described in Section \ref{sec:algorithms}. We will start from small
examples, then derive recursive formulas and finally obtain
estimations for the algorithms. Note that these derivations do not
take into account the sparsity of the matrices but are done for dense
ones. In practice, the matrices are sparse, and therefore the critical
path may be shorter.

\subsection{Critical path length for recursive inverse Cholesky}\label{cpl_estimation_rinch} 

Here we consider the critical path length of the RINCH algorithm. In Algorithm~\ref{RINCH_presudocode}, 
on every level of recursion,
there are 3 matrix-matrix multiplication operations (lines \ref{algline:rinch_line6},\ref{algline:rinch_line7} and \ref{algline:rinch_line9})
that cannot be performed in parallel and
3 matrix-scalar multiplication or matrix addition operations
(lines \ref{algline:rinch_line7} and \ref{algline:rinch_line9}).  Finally, one has to do 2 recursive calls that also cannot
be done in parallel because of data dependencies. The critical path
will for large systems be dominated by the matrix-matrix
multiplications and the recursive calls.

Let us denote the total critical path length by $\Psi(N)$, where $N$
is the matrix size. We assume that the critical path length of a
matrix-matrix multiplication task depends on the matrix size and we
will denote it as $\xi(N)$. On the lowest level of recursion, where
the leaf matrices are reached, instead of calling the task itself, one
calls leaf-level routines. We will treat these routines as tasks with
critical path length~1.

For simplicity we also assume that the leaf level matrices are of size
$1$ and that the number of non-leaf levels in the hierarchy is equal
to $L = \log_2 \left( N \right)$. Also, this leads to $\xi(1) = 1$ and
$\Psi(1) = 1$, or, in other words, all leaf-level routines
have a critical path of length 1.  Operations like matrix addition or
matrix-scalar multiplication, where each node in the matrix quad-tree
is visited only once, have critical path length $\log_2(N)+1$.  For
simplicity we will assume that the matrix dimension is equal to a
power of $2$.

For example, if we have a $2 \times 2$ matrix, then the total length
of the critical path can be written as
\begin{equation}
\Psi(2) = 3 \xi (1) + 3\cdot 1 + 2 \Psi(1) = 3 + 3 + 2 = 8.
\end{equation} If one increases the matrix size by a factor 2, then 
\begin{equation}
\Psi(4) = 3 \xi(2) + 3\left(\log_2(2) + 1\right) + 2 \Psi(2) = 3 \xi(2) + 6 + 2 \cdot 8.
\end{equation}

Let us denote by $P(N)$ the critical path length of all non-RINCH
tasks on the current level of hierarchy. In our case, $$ P(N) = 3 \xi
(N/2) + 3\left(\log_2(N/2) + 1\right), $$ i.e.\ we have 3 matrix-matrix
multiplications of $N/2 \times N/2$ matrices and 3 operations with
critical path length $\log_2(N/2) + 1$.  By $Q$ we will denote the
number of recursive RINCH calls on the current level. Note that none
of them can be done in parallel and therefore the critical path length
can then be written as
\begin{equation}
\Psi(N) = P(N) + Q \Psi(N/2), \label{rinch_recursive_psi}
\end{equation}
which gives
\begin{equation}
\begin{split}
 \Psi(N) = \sum\limits_{i = 0}^{L-1} Q^{i} P\left( \frac{N}{2^{i}} \right) + Q^L = \sum\limits_{i=0}^{L-1} 2^i P\left(\frac{N}{2^i}\right) + 2^L  \\
 = 2^{L} \sum\limits_{i=0}^{L-1} 2^{i-L} P\left( \frac{N}{2^i} \right) + 2^L.
\end{split} \label{rinch_cpl_sum}
\end{equation}

The critical path length of matrix-matrix multiplication with a
quad-tree matrix representation scales as a second order polynomial in
$\log_2(N)$ \cite{rubensson2016locality} and therefore $P(N) = c_1
\log_2^2(N) + c_2 \log_2(N) + c_3$ for some constants $c_1, c_2, c_3$
independent of $N$. Combining this with \eqref{rinch_cpl_sum} gives the critical path length for the RINCH algorithm:
\begin{equation}
\begin{split}
\Psi(N) = \left(1 + 6c_1 + 2c_2 + c_3 \right)N - c_1 \log_2^2(N) - (4c_1 + c_2)\log_2(N) \\ - 6c_1 - 2c_2 - c_3.
\end{split}
\end{equation}
For full derivation see Appendix \ref{appendixA}.  Clearly, the
critical path length grows linearly with system size, i.e.\ $\Psi(N) =
\Theta(N)$ in the sense of algorithm complexity analysis.

\subsection{Critical path length for iterative refinement} \label{cpl_estimation_irsi} Here we consider the critical path length of the IRSI algorithm.
We have already mentioned that the critical path length of
matrix-matrix multiplication with a quad-tree matrix representation
scales as a second order polynomial in $\log_2(N)$
\cite{rubensson2016locality}. Since the iterative refinement
procedure has matrix-matrix multiplication as dominating operation,
and assuming that the number of iterations does not depend on $N$, its
critical path length scales the same way.

\subsection{Critical path length for localized inverse factorization} \label{cpl_estimation_rif}

Here we derive the critical path length of the LIF algorithm.
Let us first consider a single iterative refinement iteration
(lines~\ref{algline:iterrefine_start} to \ref{algline:iterrefine_end}
of Algorithm~\ref{alg:localized_inverse_factorization}). Each
construction of the correction matrix $M_i$ and update of the inverse
factor (lines~\ref{algline:iterrefine_start}
and~\ref{algline:Z_update}) requires $m$ matrix-matrix multiplications,
$m$ matrix-scalar multiplications and $m$ matrix additions. Thus, this
update has a critical path of length $m \xi(N) + 2m(\log_2(N)+1)$,
where $\xi(N)$ is the critical path length of the matrix-matrix
multiplication and $N$ is the matrix size at the current level of the
hierarchy. The update of the $\delta$-matrix
(line~\ref{algline:delta_update}) has critical path length $2 \xi(N) +
3(\log_2(N) + 1),$ since only two matrix-matrix multiplications cannot
be done in parallel. Since $\delta$ cannot be updated until $Z$ has
been updated, the total critical path length of a single iteration is
\begin{equation}
\Psi_{it} = (m+2)\xi(N) + (2m+3)(\log_2(N)+1).
\end{equation}

It was shown in
\cite{localized_inverse_factorization} that the number of iterative
refinement iterations needed to reach an accuracy $\|\delta_{i}\|_2 <
\varepsilon$ is bounded by
\begin{equation}\label{eq:iter_refine_no_of_iters}
  K_{max} = \left\lceil\frac{\log{\left(\frac{\log \varepsilon}{\log \left(1-{\lambda_{\mathrm{min}}(S)}/{\lambda_{\mathrm{max}}(S)}\right)}\right)}}{\log{(m+1)}}\right\rceil
\end{equation}
at any level in the hierarchy.
If the condition number
${\lambda_{\mathrm{max}}(S)}/{\lambda_{\mathrm{min}}(S)}$ of $S$ is
  bounded by a constant independent of $N$ then $K_{max}$ is also
  bounded by a constant independent of $N$.
For simplicity of exposition we will in the following assume that the
number of iterative refinement iterations is exactly equal to
$K_{max}$ at all levels and branches of the recursion.  Then, the
critical path length of the localized iterative refinement procedure
is given by
\begin{equation}
\Psi_{ref}(N) = K_{max} \left( (m+2)\xi(N)  + (2m+3)(\log_2(N)+1) \right).
\end{equation}

Apart from refinement of the inverse factor, one has to prepare
initial guesses for the refinement. The construction of $\delta_0$
(line \ref{algline:delta_zero_construction} of Algorithm \ref{alg:localized_inverse_factorization}) requires 2 matrix-matrix multiplications of matrices of size
$N/2$, a separate transpose of size $N/2$ and a single multiplication
by a scalar. Therefore its critical path length is $2\xi(N/2) +
\log_2(N/2) + 1 + \log_2(N) + 1 = 2\xi(N/2) + 2 \log_2(N) + 1.$ Then,
in order to be able to prepare initial guesses, one also needs to do 2
recursive calls, each with critical path length of $\Psi(N/2).$

Similarly to the previous algorithm, let us denote by $R(N)$ the critical path length of all non-LIF tasks on the current level of the hierarchy:

\begin{equation}
R(N) = K_{max} \left( (m+2)\xi(N)  + (2m+3)(\log_2(N)+1) \right) + 2\xi(N/2) + 2 \log_2(N) + 1.
\end{equation} Then we arrive at the following recursive relation:

\begin{equation}
\Psi(N) = R(N) + Q \Psi(N/2),
\end{equation} where, similarly to the case of RINCH, $Q$ is the number of recursive calls which cannot be handled in parallel. This form is very similar to \eqref{rinch_recursive_psi} and gives us the following non-recursive relation when expanding for arbitrary $L$:

\begin{equation}
\Psi(N) = \sum\limits_{i = 0}^{L-1} Q^i R\left( \frac{N}{2^i}\right) + Q^L.
\end{equation}

Taking into account the fact that recursive LIF calls do not depend on
each other and thus can be performed in parallel, one can set $Q = 1$
and the expression becomes

\begin{equation}
\Psi(N) = \sum\limits_{i=0}^{L-1} R\left(  \frac{N}{2^i}\right)+ 1. \label{lif_cpl_sum}
\end{equation}

Using again the fact that the critical path length of matrix-matrix
multiplication with a quad-tree matrix representation scales as a
second order polynomial in $\log_2(N)$ we have that $R(N) =
c_1\log_2^2(N) + c_2\log_2(N) + c_3$ for some constants $c_1,c_2,c_3$
independent of $N$. Combining this with \eqref{lif_cpl_sum} gives the critical path length for the localized inverse
factorization algorithm:

\begin{equation}
\begin{split}
\Psi(N) = \frac{1}{2}c_1 \log_2^3(N) + \left(\frac{1}{2}c_1 + \frac{1}{2}c_2 \right)\log_2^2(N) \\ + \left(\frac{1}{6}c_1 + \frac{1}{2}c_2 + c_3\right)\log_2(N) + 1.
\end{split}
\end{equation}
For full derivations see Appendix \ref{appendixB}. One can observe
that the leading term is $\log_2^3(N)$ and $\Psi(N) =
\Theta(\log^3(N))$ in the sense of algorithm complexity analysis.  The
critical path length grows slower than linearly.

\section{Chunks and Tasks implementation} \label{sec:cht_implementation}
We implemented the RINCH, IRSI and LIF algorithms using the Chunks and
Tasks programming model~\cite{CHT-PARCO-2014}.
Chunks and Tasks is a task-based model suitable for dynamic and
hierarchical algorithms that allows a programmer to write parallel
code without thinking about synchronization or message passing between
different processes. When using this model, the
programmer writes the code by defining objects of two classes. Objects
of the first class represent pieces of data, while objects of the
second class represent pieces of work which is to be done. They are
referred to as chunks and tasks, respectively. Each task has one or more
input chunks and a single output chunk. Once a chunk has been registered to
a runtime library, it is read-only. Identifiers of chunks and tasks
are provided by the library and the user cannot modify them. This
makes the model free from deadlocks and race conditions.  In the
present work we use the Chunks and Tasks model as implemented in the
proof-of-concept runtime library CHT-MPI, freely available for downloading at \url{http://chunks-and-tasks.org/}. The library is written in C++ and
parallelized using MPI and Pthreads~\cite{CHT-PARCO-2014}. A work stealing approach is used to achieve load balancing.

The present work makes use of and extends CHTML~\cite{rubensson2016locality}
implemented using the Chunks and Tasks model.
CHTML includes a locality-aware sparse matrix-matrix
multiplication implementation that takes advantage of data locality to
reduce communication costs. 
In CHTML, a matrix is represented as a
quad-tree and all non-leaf nodes are chunks containing identifiers of
child nodes. If a submatrix is zero, then the corresponding identifier
is set to the null chunk identifier. Only on the leaf level, chunks
contain matrix elements. This approach allows to exploit sparsity
patterns dynamically by skipping operations on zero submatrices.
We use a block-sparse leaf matrix library for the representation of
submatrices at leaf level, as in~\cite{rubensson2016locality}. This
means that each leaf matrix is divided into small blocks and only
non-zero blocks are stored in memory.

The Chunks and Tasks programming model enforces certain restrictions
on how a program can be implemented.  Those restrictions exist both to
help the programmer in the development of his program and to make
the implementation of efficient runtime libraries feasible.  The
Chunks and Tasks model is purely functional and does not support
global variables. The only way a task can access data, besides compile
time constants, is through the input chunks.

Another restriction regards the uninterrupted execution of tasks. A
task cannot synchronize with child tasks during its execution.  This
means that the runtime library can rely on uninterrupted
execution of tasks which is important for efficient allocation of
resources. For the application programmer this means that while-loop
constructs where the termination condition depends on results of child
tasks registered in the while loop body cannot be implemented in the
usual way. When the number of iterations is limited an alternative is
to register a predetermined number of tasks in a for-loop including
tasks that check the termination condition. When the termination
condition has been reached, the superfluous remaining tasks
immediately return from their task executions.  Otherwise one may
encapsulate the body of the while-loop into a task and use recursive
registration of tasks.  The updated variables are encapsulated in
chunks and the chunk identifiers are provided as input to a task of
the same type together with all necessary parameters for the recursive
task registration. The next modification of the variables happens when a new task is
executed.  We used this approach for the iterative refinement procedure
as applied in the localized inverse factorization, see Algorithm \ref{IterRefine_CHT_pseudocode}. Note that this
restriction applies only to task execution. In the main program the
execution of a task is a blocking call that waits until all descendant
tasks have been executed. Thus, an algorithm like IRSI can be
implemented in the standard way.

To illustrate how the Chunks and Tasks code looks and help the reader to understand the essence of the task-based parallel programming, we include pseudocodes, which are simplified versions of the real code but include some of its key characteristics. Algorithm \ref{RINCH_CHT_pseudocode} represents  RINCH, Algorithm \ref{IterRefine_CHT_pseudocode} represents IRSI and Algorithm \ref{LIF_CHT_pseudocode} represents LIF correspondingly. The value \emph{realmax()}  in Algorithm \ref{LIF_CHT_pseudocode} represents the largest real value available, it is assigned to the error as the starting value and the error is then gradually reduced. The simplification concerns evaluations of mathematical expressions. In practice, every operation with a matrix is a task which typically generates a whole hierarchy of tasks. We write "evaluate(\emph{expr})" to indicate the registration of tasks needed to evaluate an expression \emph{expr}. Note that the variables are actually the identifiers of the corresponding chunks. As one can see, the codes consist almost purely of registrations of tasks and, in the end of the day, all tasks will be turned into chunks. The pseudocodes also ignore the necessity to derive tasks from a base task class and implement certain obligatory member functions. 

\begin{algorithm}[t]
\begin{algorithmic}[1]
\REQUIRE $S$
\ENSURE $Z$
\IF{"lowest level reached"}
	\STATE{$Z = $ InvChol($S$);}
	\STATE{\textbf{return} $Z$;}
\ELSE
	\STATE{$Z_{00} = $ registerTask(RINCH($S_{00}$));}
	\STATE{$R = $ registerTask(evaluate($Z_{00}^* S_{01}$));}
	\STATE{$Q = $ registerTask(evaluate($-R^* R + S_{11}$));}
	\STATE{$Z_{11} = $ registerTask(RINCH($Q$));}
	\STATE{$Z_{01} = $ registerTask(evaluate($-Z_{00} R Z_{11}$));}
	\STATE{$Z_{10} = NULL$;}
	\STATE{$Z = \begin{bmatrix} Z_{00} & Z_{01} \\ Z_{10} & Z_{11} \end{bmatrix}$}
	\STATE{\textbf{return} $Z$;}
\ENDIF
\end{algorithmic}
\caption{RINCH Chunks and Tasks pseudocode}
\label{RINCH_CHT_pseudocode}
\end{algorithm}

\begin{algorithm}[t]
\begin{algorithmic}[1]
\REQUIRE $Z_{prev}, S, \delta, e, e_{prev}, m$
\ENSURE $Z$
\IF{$e > e^{m+1}_{prev}$}
	\STATE{$Z = Z_{prev};$}
	\STATE{\textbf{return} $Z$;}
\ELSE
	\STATE{$M = $ registerTask(evaluate($Z_{prev} \sum\limits_{k=1}^{m}b_k \delta_{prev}^k$));}
	\STATE{$Z = $ registerTask(evaluate($Z_{prev} + M$));}
	\STATE{$\delta = $ registerTask(evaluate($\delta_{prev} - Z^* (S M) - (M^* S) Z_{prev}$));}
	\STATE{$e_{prev} = e;$}
	\STATE{$e = $ registerTask(evaluate($\|\delta\|_F$));}
	\STATE{$Z = $ registerTask(IterRefine($Z,S,\delta,e,e_{prev},m$));}
	\STATE{\textbf{return} $Z$;}
\ENDIF
\end{algorithmic}
\caption{IterRefine Chunks and Tasks pseudocode}
\label{IterRefine_CHT_pseudocode}
\end{algorithm}

\begin{algorithm}[t]
\begin{algorithmic}[1]
\REQUIRE $S,m$
\ENSURE $Z$
\IF{"switch to RINCH level reached"}
	\STATE{$Z = $ registerTask(RINCH($S$));}
	\STATE{\textbf{return} $Z$;}
\ELSE
	\STATE{$Z_{00} = $ registerTask(LIF($S_{00}$));}
	\STATE{$Z_{11} = $ registerTask(LIF($S_{11}$));}
	\STATE{$Z^{(0)} = \begin{bmatrix} Z_{00} & NULL \\ NULL & Z_{11} \end{bmatrix}$;}
	\STATE{$X = $ registerTask(evaluate($-Z_{00}^* S_{01} Z_{11}$));}
	\STATE{$\delta^{(0)} = \begin{bmatrix} NULL & X \\ X^* & NULL \end{bmatrix};$}
	\STATE{$e_{prev} = $ registerChunk(realmax());}
	\STATE{$e = $ registerChunk($e_{prev}^{m+1}$);}
	\STATE{$Z = $ registerTask(IterRefine($Z^{(0)}, S, \delta^{(0)}, e, e_{prev}, m$));}
	\STATE{\textbf{return} $Z$;}
\ENDIF
\end{algorithmic}
\caption{LIF Chunks and Tasks pseudocode}
\label{LIF_CHT_pseudocode}
\end{algorithm}

To gain more from modern C++ and make the code flexible when it comes to the choice of the algorithm version, we employed template specialization techniques for specifying the choice instead of making several similar codes. For example, the recursive inverse factorization algorithm can work with either local or non-local refinement.

\section{Experimental results}\label{sec:experiments}
To numerically investigate the scaling behaviour of the algorithms, we
used basis set overlap matrices coming from discretizations using
Gaussian basis sets for two kinds of molecules, quasi-linear Glutamic
Acid-Alanine (Glu-Ala) helices and three-dimensional water
clusters. The xyz coordinates, which were partially used
in~\cite{Ergo-JCTC-2011}, can be downloaded from
\url{http://ergoscf.org}.  The overlap matrices were generated using
the Ergo open-source program for linear-scaling electronic structure
calculations~\cite{Ergo-SoftwareX-2018}, publicly available at
\url{http://ergoscf.org} under the GNU Public License (GPL)~v3.  The
computations were done using the standard Gaussian basis set STO-3G,
with block-sparse leaf level matrix representation. The calculations
were performed in double precision, with leaf matrix size $4096$ and
blocksize $32$. The chunk cache size was set to 4 gigabytes. Moreover,
in order to enforce sparsity, truncation of small elements was
performed after every matrix-matrix multiplication with threshold
$10^{-5}$ and after all calls to leaf level routines. The truncation
procedure performs truncation of leaf internal blocks so that the
blocks with Frobenius norm smaller than the threshold value are
removed. The basis overlap matrix $S$ was also truncated with the same
threshold, and then factorized.  We have already mentioned that in the
localized inverse factorization, one may, instead of going down to the
leaf level, employ RINCH somewhere higher up in the hierarchy.  We did
so when the current matrix size became smaller than or equal to
16384. For the iterative refinement we used $m=4$. 
The local version of the refinement \eqref{refinement_local} was used in
all calculations. When performing iterative refinement with scaled
identity as starting guess (IRSI), we did not take the time needed to
compute the Gershgorin bounds into account.

All the computations were performed on the Beskow cluster located at the
PDC high performance computing center at the KTH Royal Institute of
Technology in Stockholm. The test calculations were performed using a
development version of CHT-MPI compiled with the GCC 7.3.0 g++
compiler, Cray MPICH 7.7.0 and OpenBLAS 0.2.20 \cite{OpenBLAS} for leaf-level
routines. The OpenBLAS was configured to use a single thread. The
Beskow cluster consists of 2060 nodes, where each node has 2 Intel
Xeon E5-2698v3 CPUs with 16 cores each running at 2.3 GHz frequency
(3.6 GHz with Turbo-boost), with 64 GB of RAM. The nodes are connected
with Cray Aries high speed network of Dragonfly topology. In the
CHT-MPI library every process uses a set of threads of which some
execute tasks and some do auxiliary work like work stealing. We used a
single process per node, 31 worker threads per process thus leaving a
single core for communication and other auxiliary work.

In some of the figures presented in the following subsections, some data
points are missing for one or several methods. This indicates that the
calculation was terminated due to insufficient memory. We noted that
this happened more often for the RINCH and IRSI methods.

\subsection{Scaling with system size} 

\begin{figure}[h!]
\begin{minipage}[h]{0.49\linewidth}
\center{\includegraphics[width=0.9\linewidth]{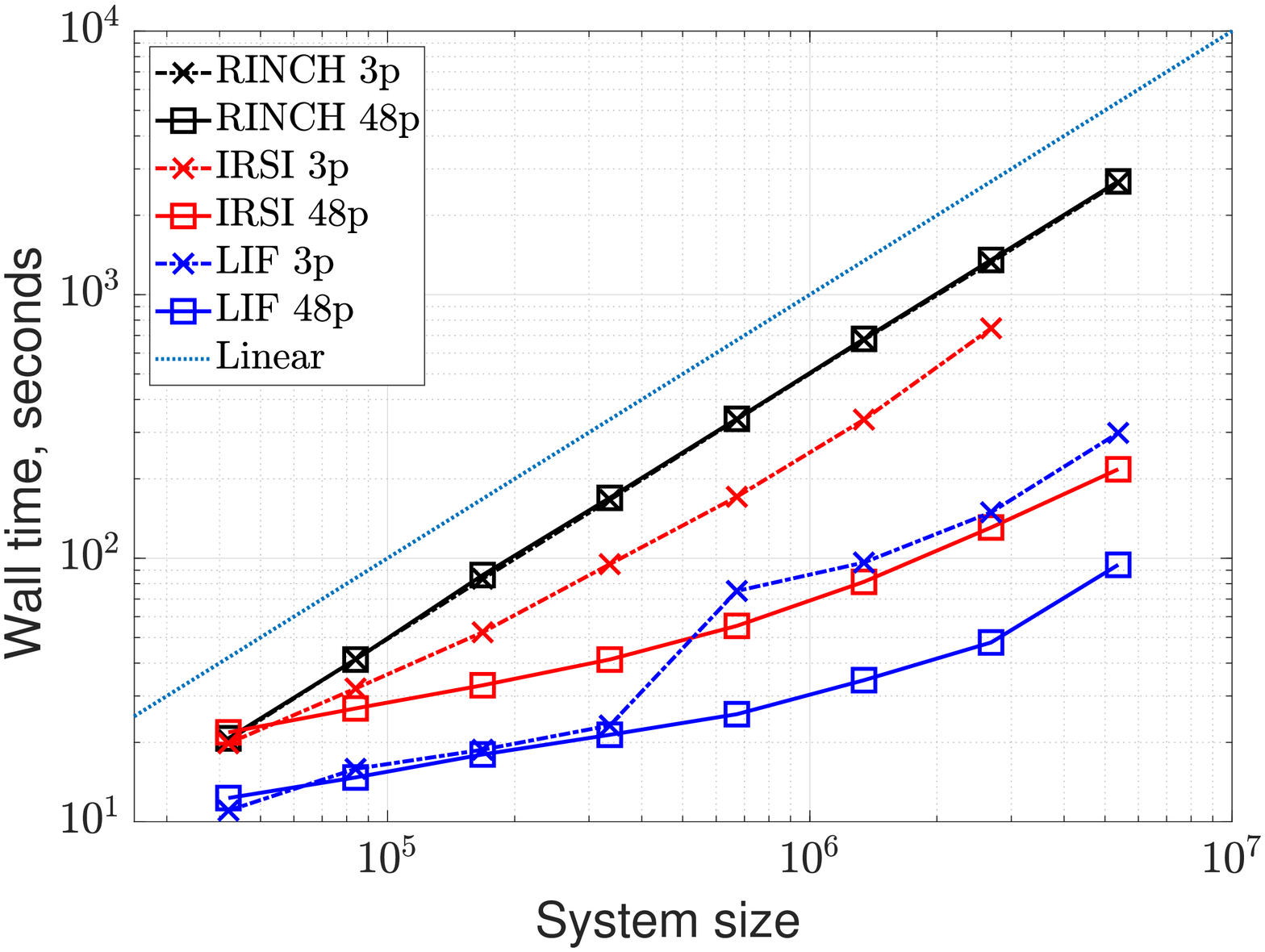} \\}
\end{minipage}
\hfill
\begin{minipage}[h]{0.49\linewidth}
\center{\includegraphics[width=0.9\linewidth]{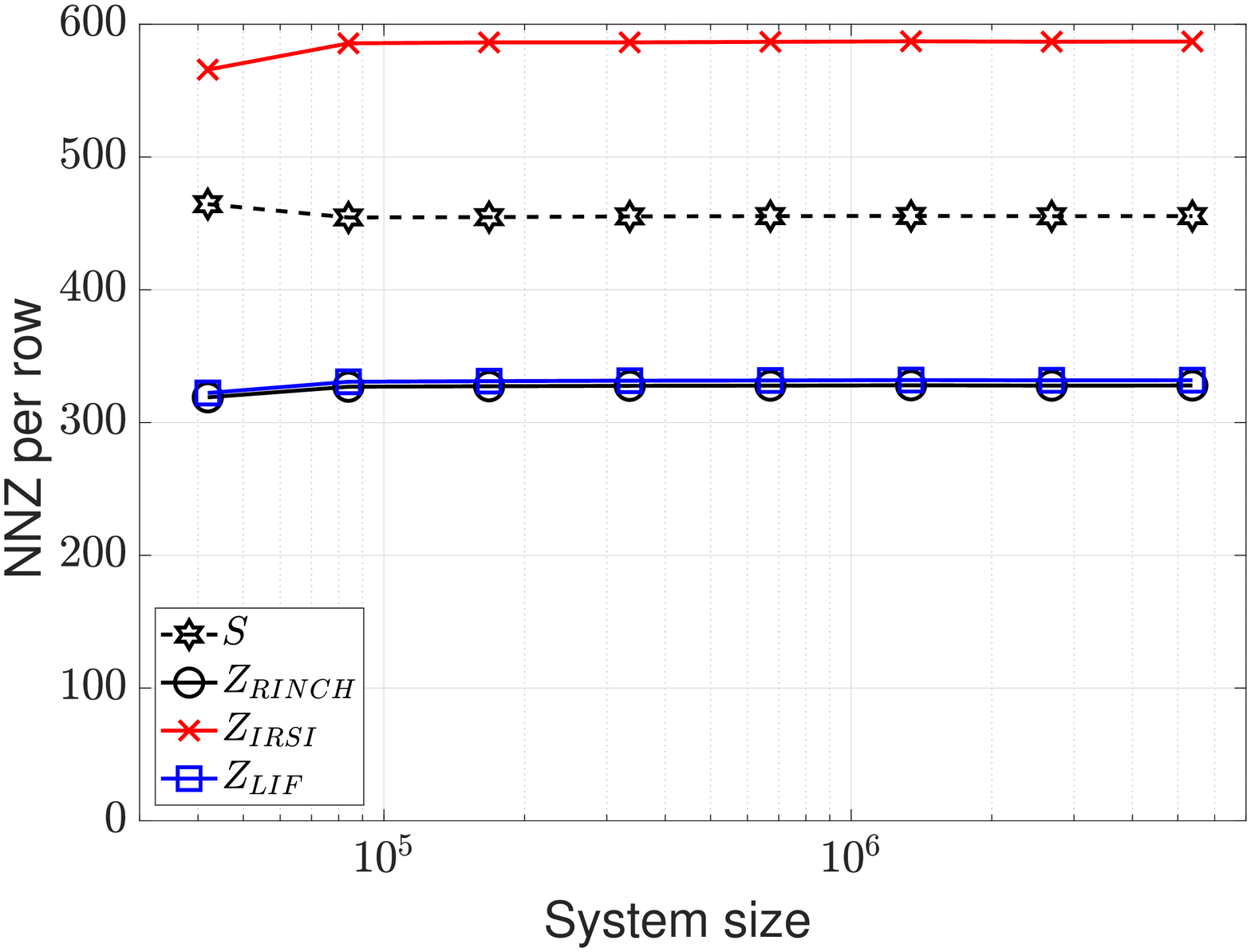} \\}
\end{minipage}
\caption{Left panel: Scaling with system size of the RINCH, IRSI and LIF algorithms for Glu-Ala helices of increasing length. The tests were made for 3 and 48 processes involved. Right panel: numbers of non-zero elements per row in the corresponding inverse factors and in the original overlap matrix $S.$ }
\label{ris:image1}
\end{figure} 

Figure \ref{ris:image1} demonstrates how the algorithms behave with
increasing system size in case of quasi-linear Glu-Ala helices. The
curves represent scaling with 3 and 48 processes used. One can see
that RINCH scales linearly, while LIF gives better time and scales
better than linearly. IRSI shows similar scaling results, but works
slower. One possible explanation of the scaling of LIF and IRSI is
that before being saturated with work, processes have extra resources
available. When saturated, they demonstrate a behaviour of the same
type as RINCH shows. The number of non-zeros (NNZ) per row in $S$
stays almost constant, as well as in inverse factors computed with all
three algorithms. IRSI provides the densest factors, while LIF and
RINCH give approximately the same sparsity. Note that RINCH gains
almost nothing from increasing the number of processes.

\begin{figure}[h!]
\begin{minipage}[h]{0.49\linewidth}
\center{\includegraphics[width=0.9\linewidth]{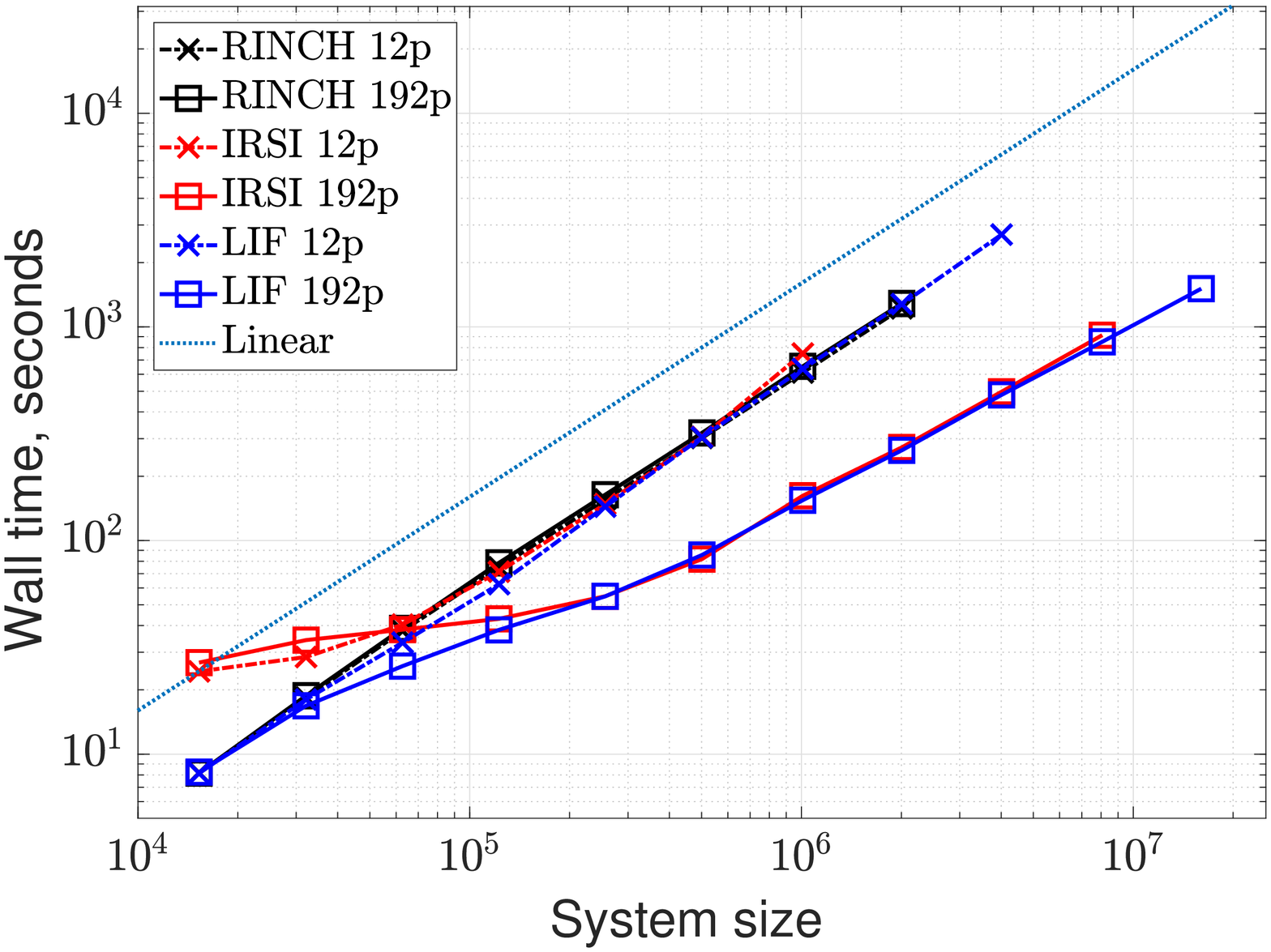} \\}
\end{minipage}
\hfill
\begin{minipage}[h]{0.49\linewidth}
\center{\includegraphics[width=0.9\linewidth]{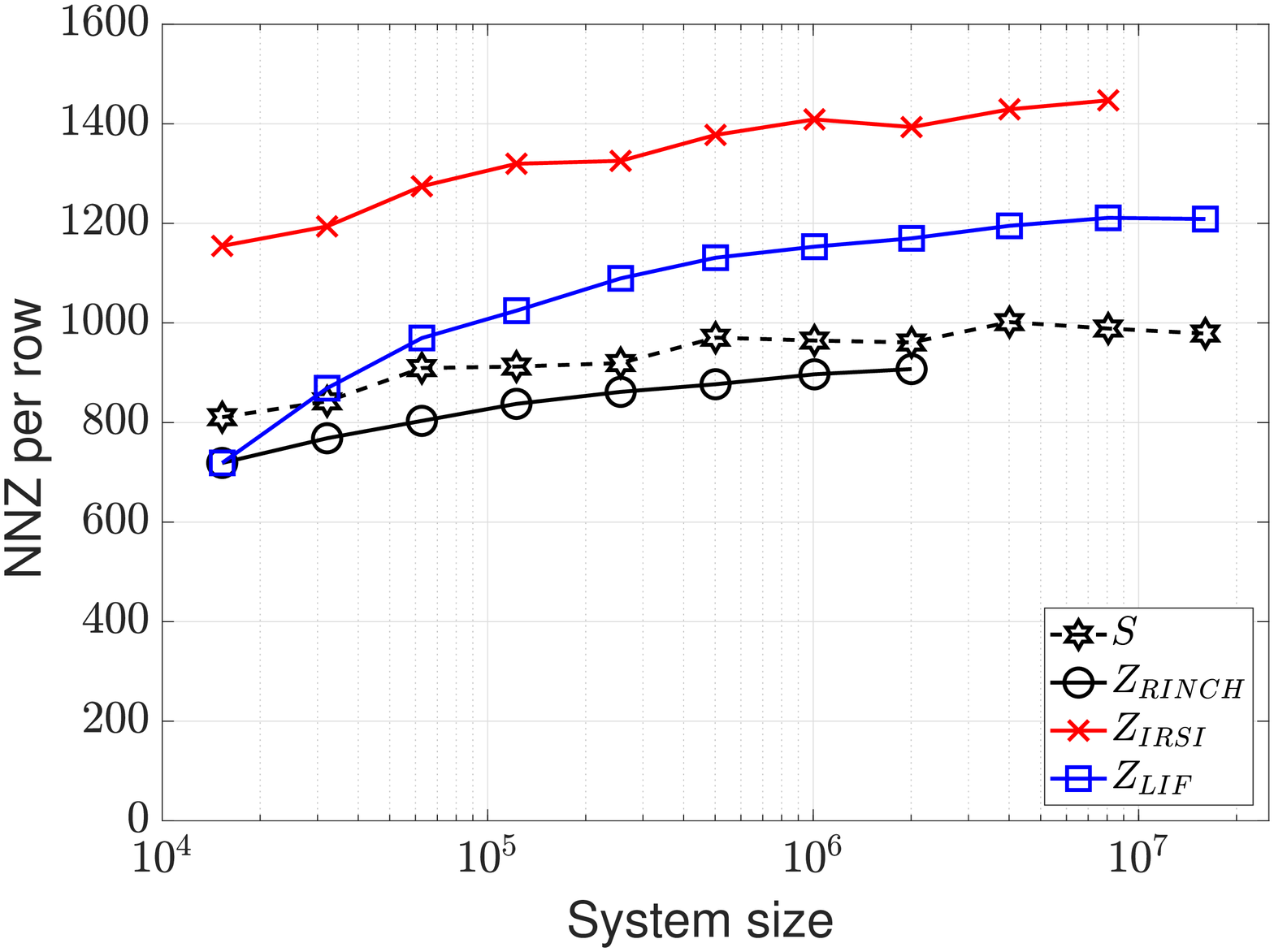} \\}
\end{minipage}
\caption{Left panel: Scaling with system size of the RINCH, IRSI and LIF algorithms for water clusters of increasing size. The tests were made for 12 and 192 processes involved. Right panel: numbers of non-zero elements per row in the corresponding inverse factors and in the original overlap matrix $S.$ }
\label{ris:image2}
\end{figure}

Figure \ref{ris:image2} demonstrates how the algorithms behave with
increasing system size in case of water clusters. Now the curves on
the left plot represent scaling with 12 and 192 processes. One can
observe linear scaling as the system size is increased. RINCH again
shows no speedup when going from 12 to 192 processes.  In case of 192
processes, LIF and IRSI scale slightly better than linearly, which
indicates that the number of processes is too large and not all of
them are saturated with work. Eventually the LIF and the IRSI
algorithms arrive at very close timings. The number of nonzeros per
row is no longer constant, it grows slowly with system size but
flattens out for large systems
\cite{localized_inverse_factorization}. Note that the RINCH algorithm was not able to handle the largest systems.
 
\subsection{Strong scaling} 
\begin{figure}[h!]
\begin{minipage}[h]{0.49\linewidth}
\center{\includegraphics[width=0.9\linewidth]{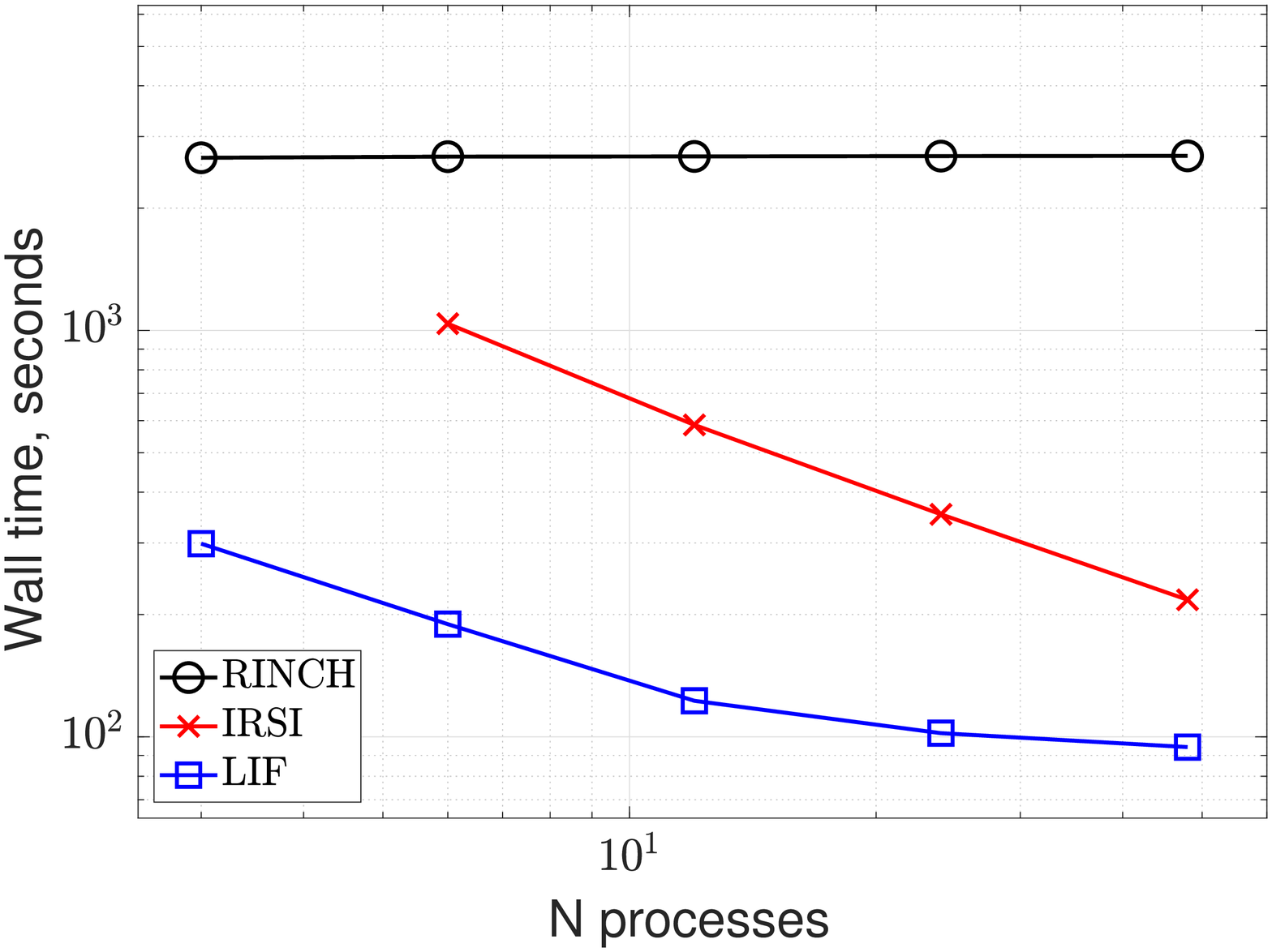} \\}
\end{minipage}
\hfill
\begin{minipage}[h]{0.49\linewidth}
\center{\includegraphics[width=0.9\linewidth]{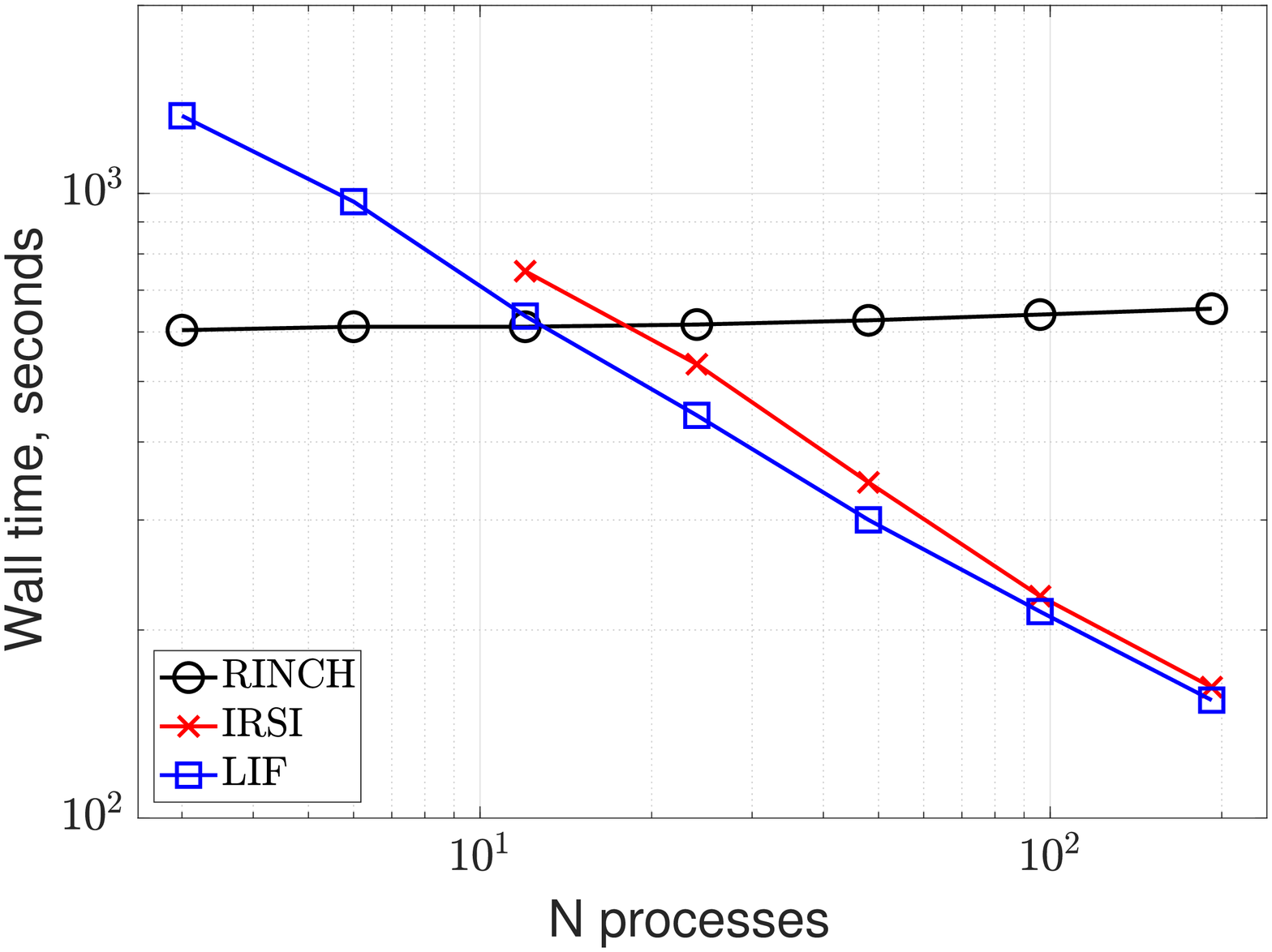} \\}
\end{minipage}
\caption{Left panel: Strong scaling of the RINCH, IRSI and LIF algorithms for a Glu-Ala helix containing 1703938 atoms, which gave a system size of 5373954 basis functions. Right panel: strong scaling of the RINCH, IRSI and LIF algorithms for a water cluster containing 432498 atoms, which gave a system size of 1009162 basis functions. For both cases, the number of processes was doubled each time while the system size was kept the same.}
\label{ris:image3}
\end{figure}

One can also fix the problem size and see how the execution time
changes when the number of processes is increased. Figure
\ref{ris:image3} shows the scaling in the strong sense for the
algorithms. As one can see, RINCH has an almost flat curve with no
speedup. The other two algorithms exhibit strong scaling with similar
behaviour, although the IRSI algorithm failed to handle small cases because of memory limitations.

\subsection{Weak scaling, critical path length, and data movement} 
Due to the task-based nature of the Chunks and Tasks model and dynamic
load balancing based on work stealing, it is very difficult to write
down the exact formula for the speedup as a function of the number of
processes involved. Therefore it is natural to use a different tool,
like a critical path estimation, in order to predict the scaling in
the weak sense. By weak scaling we mean how the parallel execution
time changes when the workload increases in direct proportion to 
the number of processes.

\begin{figure}[h!]
\begin{minipage}[h]{0.49\linewidth}
\center{\includegraphics[width=0.9\linewidth]{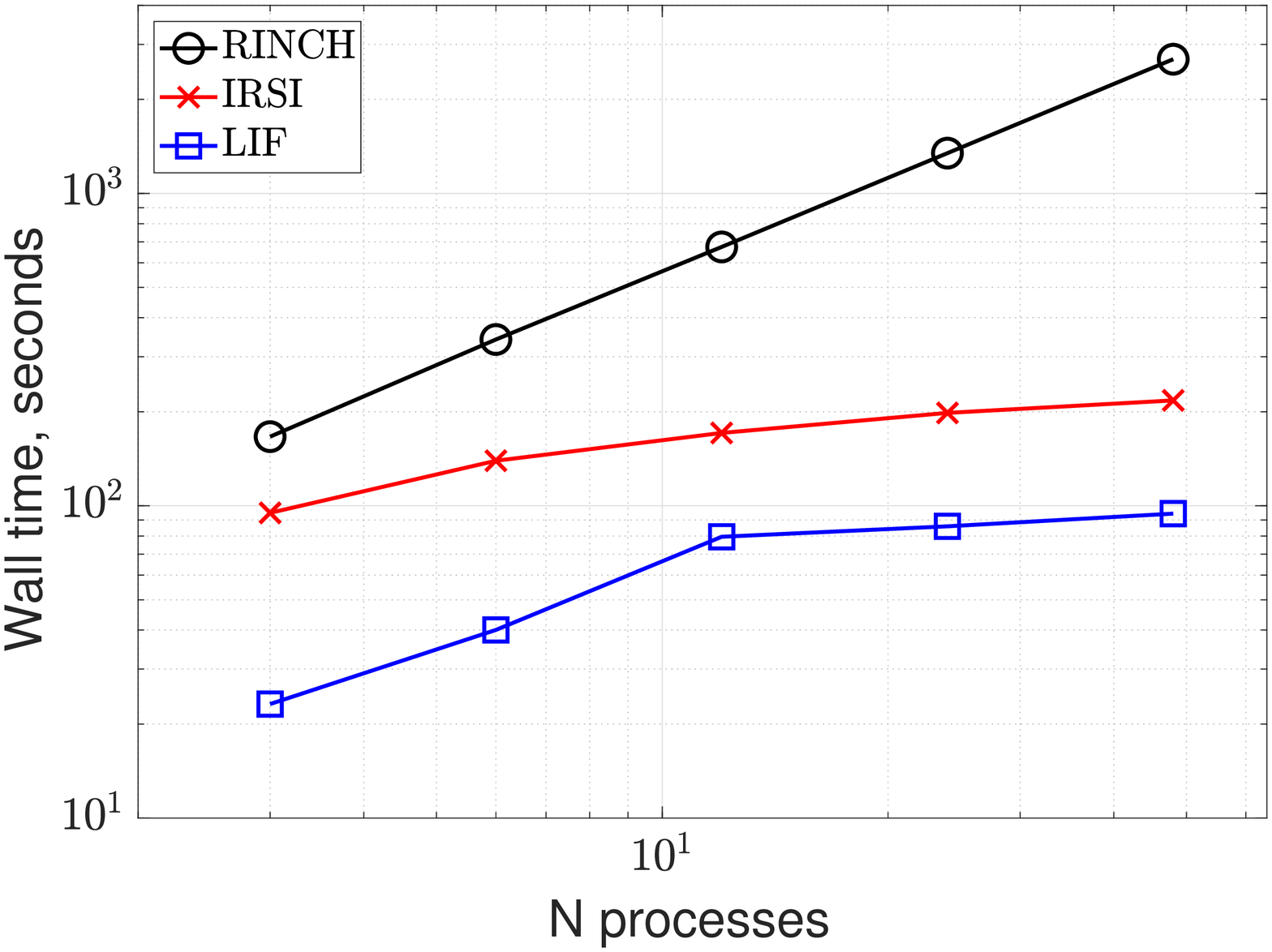} 
}
\end{minipage}
\hfill
\begin{minipage}[h]{0.49\linewidth}
\center{\includegraphics[width=0.9\linewidth]{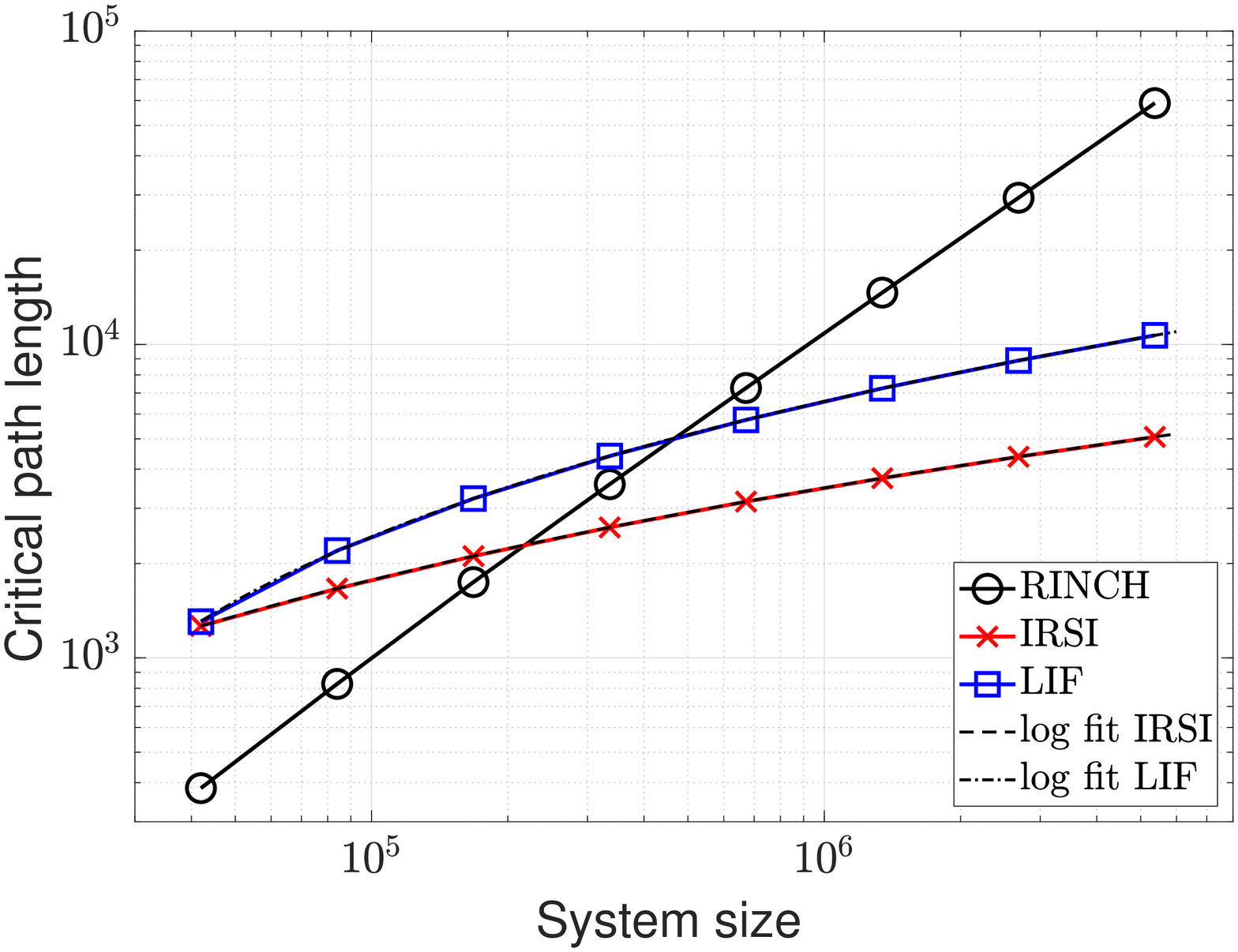} 
}
\end{minipage}
\caption{Left panel: Approximate weak scaling of the RINCH, IRSI, and
  LIF algorithms for Glu-Ala helices of increasing length. The number
  of basis functions per process was approximately fixed to $112
  \times 10^3$, so that the system size is scaled up together with the
  number of processes. Right panel: Critical path length as reported
  by the CHT-MPI library, defined as the largest number of tasks that
  have to be executed serially.  The dashed
  and dashed-dotted help lines show $c_0 + c_1 \log(N) + c_2
  \log^2(N) + c_3 \log^3(N)$ least squares fits for IRSI and LIF,
  respectively.  The data in the left panel corresponds to the 5
  rightmost points in the right plot.}
\label{ris:image4}
\end{figure} 

\begin{figure}[h!]
\begin{minipage}[h]{0.49\linewidth}
\center{\includegraphics[width=0.9\linewidth]{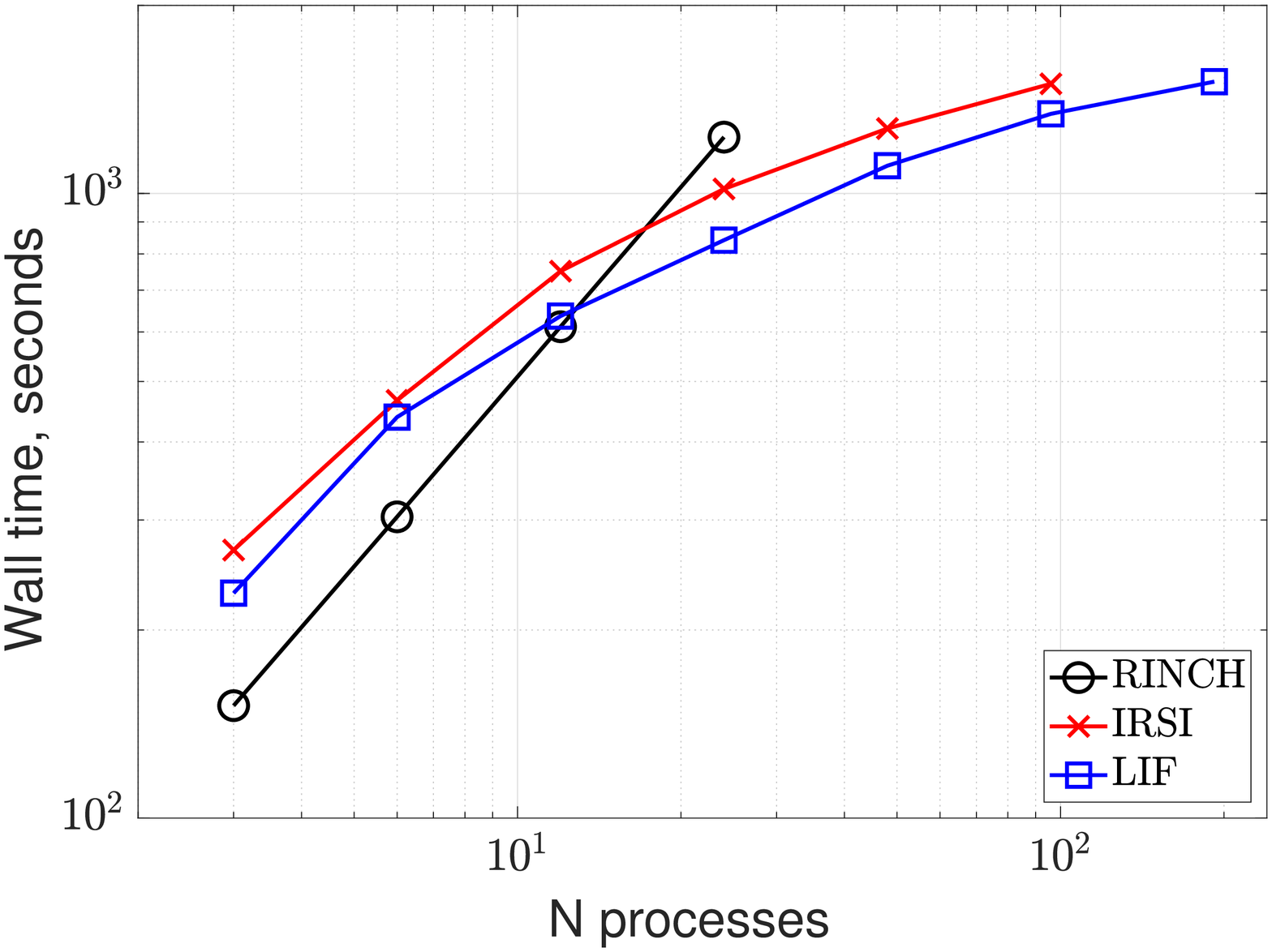} \\}
\end{minipage}
\hfill
\begin{minipage}[h]{0.49\linewidth}
\center{\includegraphics[width=0.9\linewidth]{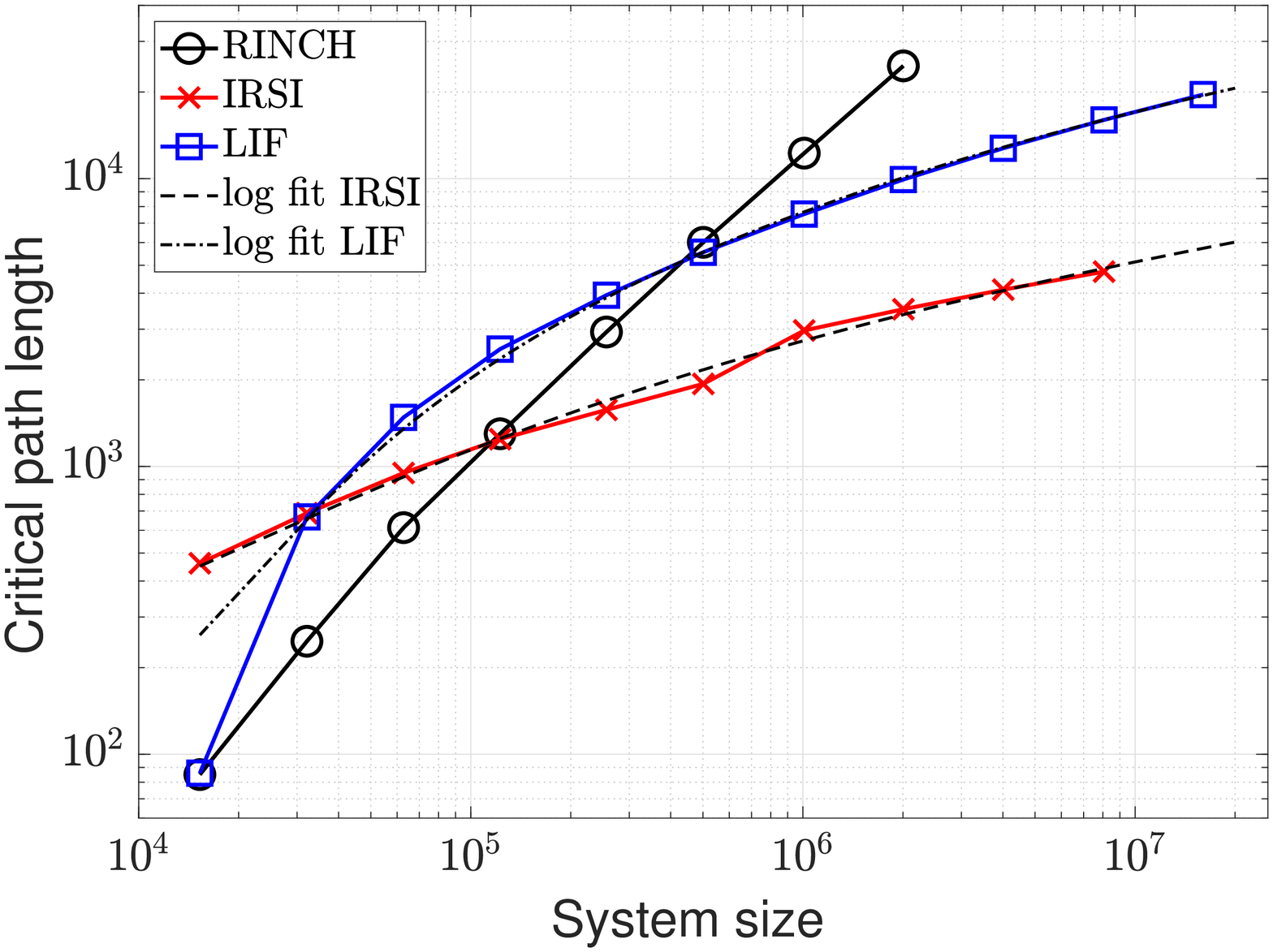} \\}
\end{minipage}
\caption{Left panel: Approximate weak scaling of the RINCH, IRSI, and
  LIF algorithms for water clusters of increasing length. The number
  of basis functions per process was approximately fixed to $84
  \times 10^3$, so that the system size is scaled up together with the
  number of processes. Right panel: Critical path length as reported
  by the CHT-MPI library, defined as the largest number of tasks that
  have to be executed serially.  The dashed
  and dashed-dotted help lines show $c_0 + c_1 \log(N) + c_2
  \log^2(N) + c_3 \log^3(N)$ least squares fits for IRSI and LIF,
  respectively.  The data in the left panel corresponds to the 7
  rightmost points in the right plot.}
\label{ris:image5}
\end{figure} 

Figures~\ref{ris:image4} and~\ref{ris:image5} show how the algorithms
scale in the weak sense and critical path lengths for Glu-Ala helices helices and water clusters,
respectively. The figures show an approximate weak scaling since,
strictly speaking, it is the number of basis functions per process and
not the workload per process that is fixed as the number of processes
is increased. Note however that the workload per process should
approach a constant in the limit of large systems.
The critical path length for the recursive inverse
Cholesky algorithm increases linearly with system size as predicted in
Section~\ref{sec:cpl_estimation}.  The critical path lengths for iterative
refinement with scaled identity and localized inverse factorization
also agree with the theoretical results saying that the critical path
lengths should increase as polynomials in $\log(N)$. In all cases the
wall times in these weak scaling tests are dominated by the execution
of tasks, including associated communication, along the critical
path.
We noted that the coefficient $c_3$ in front of the $\log^3(N)$ term
in the least squares fit was very small for both IRSI and LIF. For
IRSI this is in perfect agreement with the theory, and for LIF it can
be explained by matrix sparsity.
Although the critical paths for the localized
inverse factorization algorithm are typically longer than the critical
paths for iterative refinement with scaled identity, the localized
inverse factorization is faster.  The reason is mainly that less
communication is needed in the localized inverse factorization, as
anticipated in the introduction.

\begin{figure}[h!]
\begin{minipage}[h]{0.49\linewidth}
\center{\includegraphics[width=0.9\linewidth]{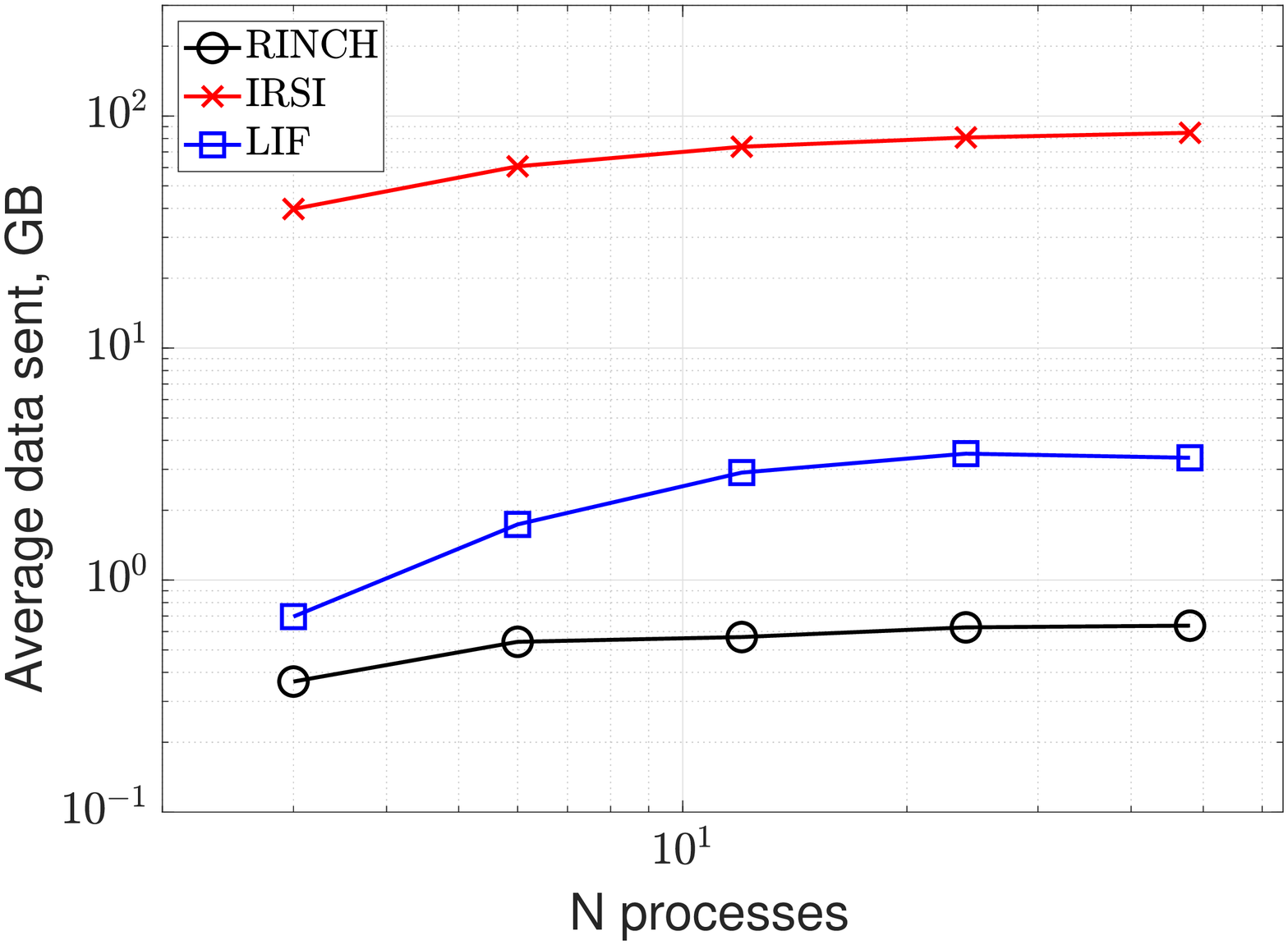} \\}
\end{minipage}
\hfill
\begin{minipage}[h]{0.49\linewidth}
\center{\includegraphics[width=0.9\linewidth]{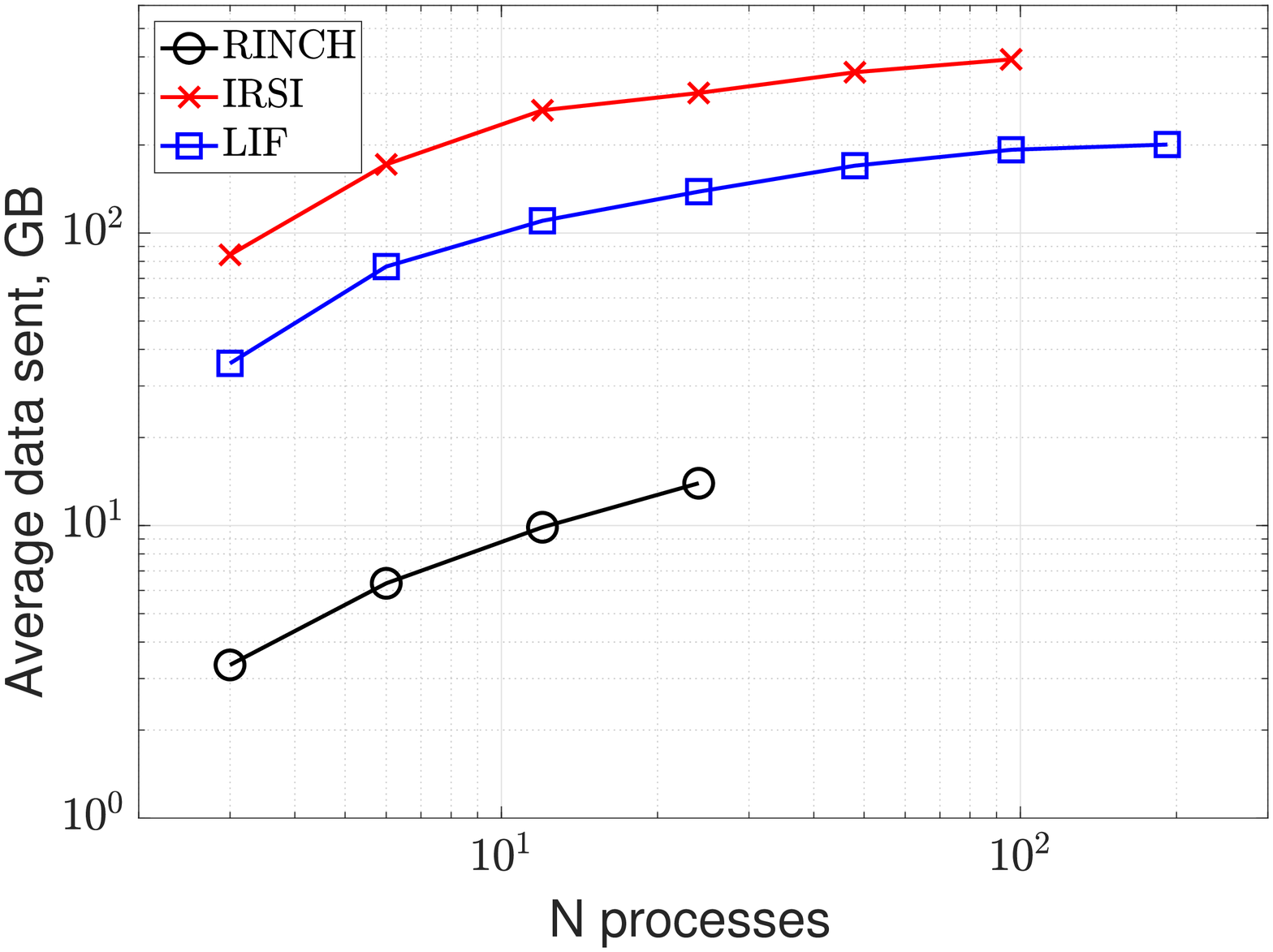} \\}
\end{minipage}
\caption{Left panel: Average amount of data sent per process for the RINCH, IRSI and LIF algorithms for the Glu-Ala helices of increasing size. The number of basis functions per process was approximately fixed to $112 \times 10^3,$ so that the system size is scaled up with the number of processes. Right panel: Average amount of data sent per process  for the RINCH, IRSI and LIF algorithms for the water clusters of increasing size. The number of basis functions per process was approximately fixed to $84 \times 10^3.$}
\label{ris:image6}
\end{figure} 

\begin{figure}[h!]
\begin{minipage}[h]{0.49\linewidth}
\center{\includegraphics[width=0.9\linewidth]{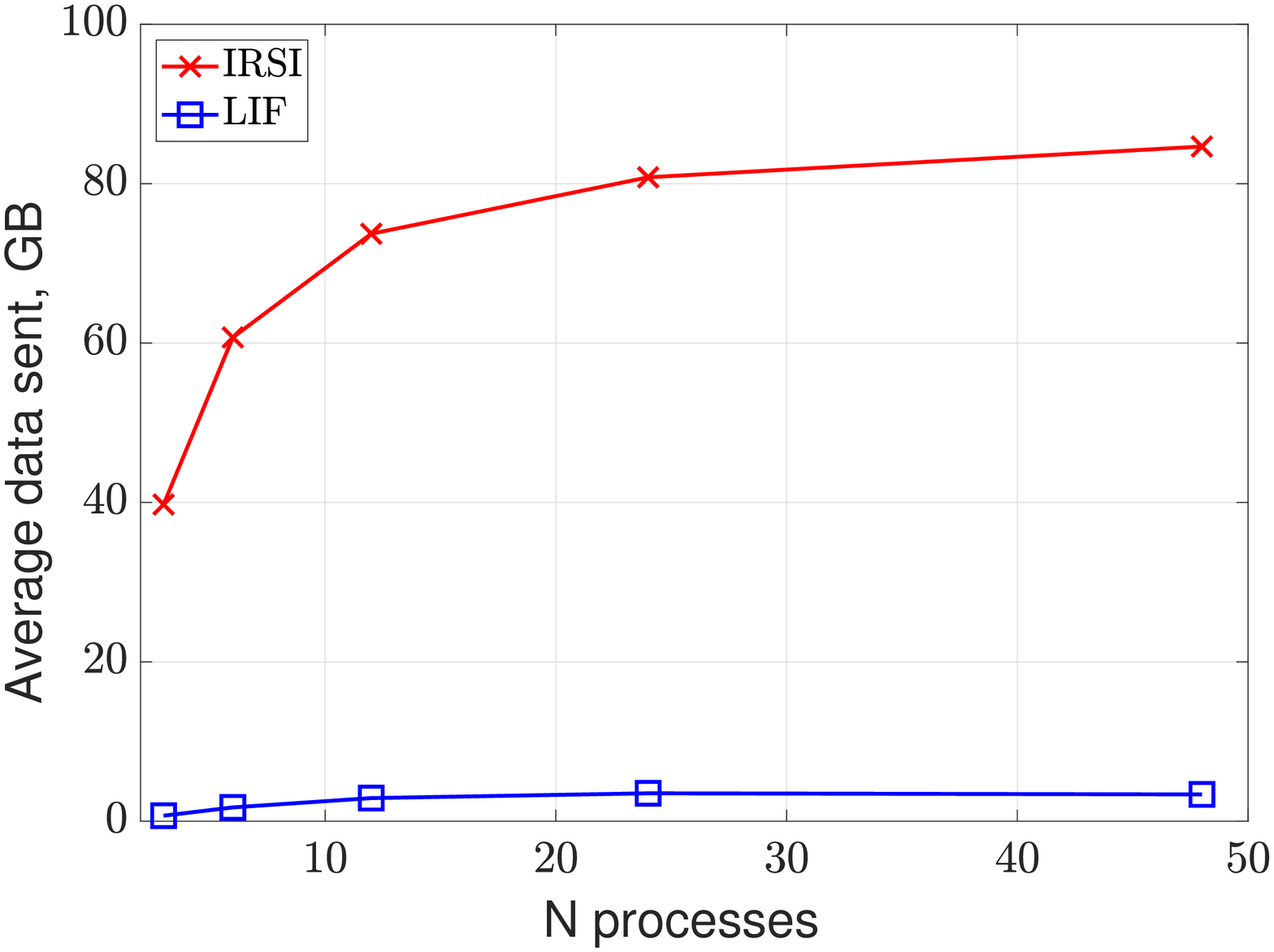} \\}
\end{minipage}
\hfill
\begin{minipage}[h]{0.49\linewidth}
\center{\includegraphics[width=0.9\linewidth]{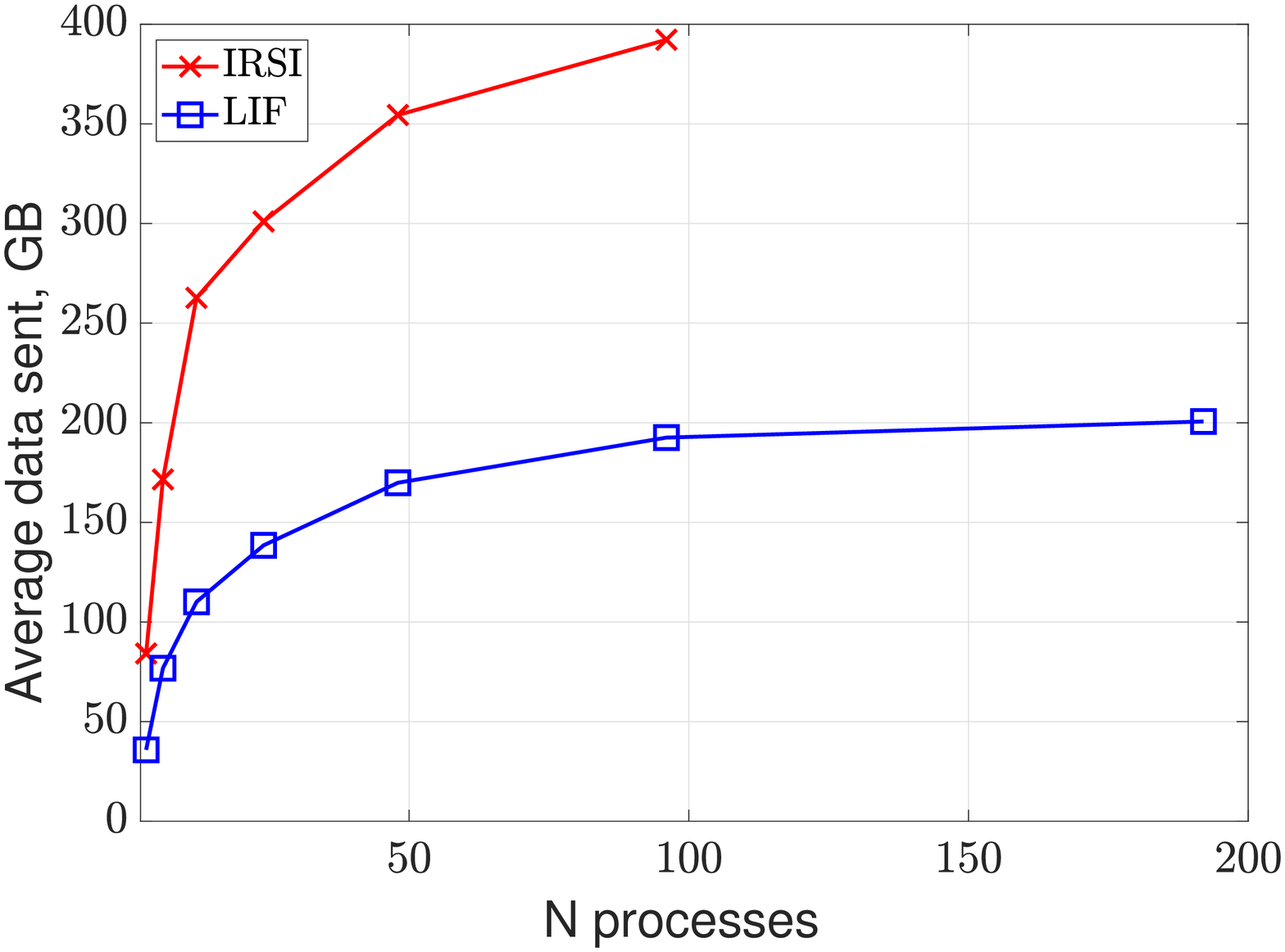} \\}
\end{minipage}
\caption{Left panel: Average amount of data sent per process for the IRSI and LIF algorithms for the Glu-Ala helices of increasing size, linear-linear scale. Right panel: Average amount of data sent per process for the IRSI and LIF algorithms for the water clusters of increasing size, linear-linear scale.}
\label{ris:image7}
\end{figure} 

Figures~\ref{ris:image6} and~\ref{ris:image7} show the average amount
of data sent per process for the weak scaling tests in
Figures~\ref{ris:image4} and~\ref{ris:image5}.
The RINCH algorithm has the smallest average amount
of data moved. This agrees with previous observations regarding the
scaling in the weak and the strong sense and the critical path growth,
since all the work is actually done by a single process only. So,
there is not so much to move. In turn, the IRSI algorithm has the
largest amount of data moved.  Localized inverse factorization
preforms better than IRSI. Although both IRSI and LIF are based on
matrix-matrix multiplication, the localized computations used in LIF
lead to a
reduced amount of communication clearly visible in the
plots.

Figure \ref{ris:image7} shows the same dependency, but in
linear-linear scale, and the advantage of LIF over IRSI is clearly
visible in both cases. Note that for the LIF algorithm, the average
amount of data moved becomes almost a constant for the system size 
increasing simultaneously with the number of processes at some point.

\subsection{Factorization error} 
Table \ref{error_table} demonstrates the difference between the
algorithms in terms of factorization error $I - Z^{*} S Z.$ We pick 2
particular cases, one for Glu-Ala helices and one for water clusters,
b.f. stands for "basis functions", i.e. the system size. One can see
that the RINCH algorithm provides the most accurate results in both
cases. In turn, IRSI gives the least accurate results. Localized
inverse factorization provides results of intermediate quality.

\begin{table}[]
\centering
\caption{Frobenius norm of the factorization error ($\|I-Z^TSZ\|_F$) for the three methods using a threshold value $10^{-5}$ for the removal of small matrix entries. See the text for details.}
\label{error_table}
\begin{tabular}{|l|l|l|l|}
\hline
                                    & RINCH                 & LIF                   & IRSI                 \\ \hline
Glu-Ala, 5373954 b.f.       		 & 0.00204 & 0.00259 & 0.02352 \\ \hline
Water clusters, 2006214 b.f.        & 0.00603 & 0.00999 & 0.02628 \\ \hline
\end{tabular}

\end{table}

\section{Discussion and concluding remarks} \label{sec:discussion}
The LIF algorithm can in principle be used
by itself all the way down to single matrix elements. However, one may
at any level in the hierarchy switch to some other algorithm to
compute an inverse factor. In our numerical examples we switched to
the recursive inverse Cholesky algorithm when the matrix size became
smaller than or equal to 16384. In this way we combined the speed of
the inverse Cholesky algorithm for small systems with the good
parallel scaling and localization properties of the LIF algorithm.
The advantage of LIF in this sense is clearly seen in Figure~\ref{ris:image2}
where LIF is the most efficient algorithm for system sizes ranging over four orders of magnitude.

A drawback of IRSI is that it requires knowledge of at least some
estimates to extremal eigenvalues to construct a starting guess. The
convergence depends on the accuracy of those estimates. The extremal
eigenvalues are not always cheap to calculate and in the present work
we used the approximate scaling from~\cite{jansik2007linear} based on
an estimate of the largest eigenvalue using the Gershgorin circle
theorem.  This approximate scaling resulted in convergence in all our
numerical experiments, besides the calculations that were killed due
to insufficient memory.  However, one should be aware that the use of
this scaling relies on an overestimation of the eigenvalue since an
exact eigenvalue in \eqref{eq:scaling_gers} would lead to divergence.

An alternative not considered in this work is the family of
multifrontal methods~\cite{duff1983multifrontal}. Similarly to the
recursive and localized inverse factorization, multifrontal methods
make use of a hierarchical decomposition of the matrix.  In the
multifrontal methods the matrix is seen as the adjacency matrix of a
graph, which is split into two disconnected subgraphs using a vertex
separator. This results in a 3 by 3 block division of the matrix.  The
2 by 2 block division of the matrix in the localized inverse
factorization can similarly be seen as the result of an edge separator
of the graph.
Since the two subgraphs in the multifrontal methods need to be
disconnected, sparsity is needed to get parallelism. Without sparsity,
there are no disconnected subgraphs. In LIF there is always
parallelism, even if there is little or no sparsity. In the
multifrontal methods explicit information about the sparsity pattern
is needed to make the split. In LIF, any split leads to convergence,
although good splits result in smaller initial factorization errors
leading to faster convergence and/or more localized computations.

Finally we would like to stress that the present work is the first
parallel implementation of the LIF algorithm and that previous
implementations of IRSI made use of Cannon's algorithm or SpSUMMA
giving an unfavorable scaling in the weak sense. Therefore, to our
best knowledge, our LIF and IRSI implementations using the Chunks and
Tasks programming model are the first parallel inverse factorization
implementations that achieve a weak scaling performance with the
computational time increasing as a low order polynomial in $\log(N)$.

\section{Acknowledgements}
Support from the Swedish national strategic e-science research program
(eSSENCE) is gratefully acknowledged. Computational resources were
provided by the Swedish National Infrastructure for Computing (SNIC)
at the PDC Center for High Performance Computing at the KTH Royal
Institute of Technology in Stockholm.

\appendix
\section{Critical path length derivation for the RINCH algorithm}\label{appendixA}

Let us first derive some useful expressions:

\begin{align}
   \sum\limits_{i=0}^{L-1}2^{i-L} ={}&  \sum\limits_{j=1}^{L} \frac{1}{2^j} = 1 - \frac{1}{2^L} = 1 - \frac{1}{N}. \label{appendix_eq1}
\\
    \sum\limits_{i=0}^{L-1}i\cdot 2^{i-L} ={}& \sum\limits_{j=1}^{L} \frac{L - j}{2^j} = L\sum\limits_{j=1}^{L} \frac{1}{2^j} - \sum\limits_{j=1}^{L}\frac{j}{2^j} \nonumber \\ 
         & = L \left(1 - \frac{1}{2^L}\right) - \left(-\frac{L}{2^L}+2-\frac{2}{2^L} \right) = \log_2(N) - 2 + \frac{2}{N}. \label{appendix_eq2}
\\
    \sum\limits_{i=0}^{L-1}i^2\cdot 2^{i-L} ={}& \sum\limits_{j=1}^{L} \frac{\left(L - j\right)^2}{2^j} = L^2 \sum\limits_{j=1}^{L} \frac{1}{2^j} + \sum\limits_{j=1}^{L} \frac{j^2}{2^j} - 2L \sum\limits_{j=1}^{L} \frac{j}{2^j} \nonumber \\
    	 & = L^2 \left(1 - \frac{1}{2^L} \right) + 6 - \left(L^2 + 4 L + 6 \right)\frac{1}{2^L} - 2L \left(-\frac{L}{2^L} + 2 - \frac{2}{2^L} \right) \nonumber \\
    	 & = L^2 + 6 - \frac{6}{2^L} - 4L = \log_{2}^{2}(N) - 4\log_2(N) + 6 - \frac{6}{N}.\label{appendix_eq3}
\end{align}

Since $\xi(N)$ is a second order polynomial in $\log_2(N)$, $P(N) = c_1 \log^{2}_2(N) + c_2 \log_2(N) + c_3.$ Then, \eqref{rinch_cpl_sum} gives, 

\begin{align}
	\Psi(N) = {}& 2^L \sum\limits_{i=0}^{L-1}2^{i-L}P\left( \frac{N}{2^i} \right) + 2^L \nonumber \\
			& = N  + N \sum\limits_{i=0}^{L-1} 2^{i-L}\left( c_1 \log_2^2 \left( \frac{N}{2^i} \right) + c_2 \log_2 \left(\frac{N}{2^i} \right) + c_3 \right) \nonumber \\
			& = N + N \sum\limits_{i=0}^{L-1} 2^{i-L}\left( c_1 \left(\log_2(N)-i \right)^2 + c_2  \left(\log_2(N)-i \right) + c_3 \right) \nonumber \\
			& = N + c_1 N \log_2^2(N)\sum\limits_{i=0}^{L-1}2^{i-L} + c_1 N \sum\limits_{i=0}^{L-1}i^2 \cdot 2^{i-L} \nonumber \\ 
		    & - 2 c_1 N \log_2(N) \sum\limits_{i=0}^{L-1} i \cdot 2^{i-L} + c_2 N \log_2(N) \sum\limits_{i=0}^{L-1} 2^{i-L} \nonumber \\
		    & - c_2 N \sum\limits_{i=0}^{L-1} i \cdot 2^{i-L} + c_3 N \sum\limits_{i=0}^{L-1} 2^{i-L} \nonumber \\
		    & = N + c_1 N \log_2^2(N) \left(1 - \frac{1}{N} \right) + c_1 N \left(\log_{2}^{2}(N) - 4\log_2(N) + 6 - \frac{6}{N}  \right) \nonumber \\
		    & - 2 c_1 N \log_2(N) \left( \log_2(N) - 2 + \frac{2}{N} \right) + c_2 N \log_2 \left(1 - \frac{1}{N} \right) \nonumber \\
		    & - c_2 N \left( \log_2(N) - 2 + \frac{2}{N} \right) + c_3 N \left( 1 - \frac{1}{N} \right) \nonumber \\
		    & = \left(1 + 6c_1 + 2c_2 + c_3 \right)N - c_1 \log_2^2(N) - (4c_1 + c_2)\log_2(N) \nonumber \\ 
		    & - 6c_1 - 2c_2 - c_3.  \label{appendix_eq4}
\end{align}

\section{Critical path length derivation for the RIF algorithm} \label{appendixB}

Auxiliary formulas:

\begin{align}
    \sum\limits_{i=1}^{L} i ={}& \frac{\left(L-1\right)L}{2}. \label{appendix_eq5}
\\
    \sum\limits_{i=0}^{L-1} i^2 ={}& \frac{L^3}{6} - \frac{L^2}{2} + \frac{L}{6}. \label{appendix_eq6}
\end{align}

Recall that $R(N)$ is a second order polynomial in $\log_2(N),$ $R(N) = c_1 \log^{2}_2(N) + c_2 \log_2(N) + c_3.$ Then, \eqref{lif_cpl_sum} gives

\begin{align}
	\Psi(N) ={}& \sum\limits_{i=0}^{L-1} R\left(\frac{N}{2^i}\right) + 1 \nonumber \\
		    & = \sum\limits_{i=0}^{K-1}\left(c_1 \log^{2}_2\left(\frac{N}{2^i}\right) + c_2 \log_2\left(\frac{N}{2^i}\right) + c_3 \right) + 1 \nonumber \\
		    & = \sum\limits_{i=0}^{L-1}\left(c_1\left(\log_2(N)-i\right)^2 + c_2 \left(\log_2(N)-i\right) + c_3 \right) + 1 \nonumber \\
		    & = c_1 \log_2^2(N)\sum\limits_{i=0}^{L-1} 1 + c_1\sum\limits_{i=0}^{L-1}i^2 - 2 c_1 \log_2(N) \sum\limits_{i=0}^{L-1} i + c_2 \log_2(N) \sum\limits_{i=0}^{L-1}1 \nonumber \\
		    & - c_2 \sum\limits_{i=0}^{L-1} i + c_3 \sum\limits_{i=0}^{L-1} 1 + 1 \nonumber \\
		    & = c_1 \log_2^3(N) + c_1 \left(\frac{L^3}{3} - \frac{L^2}{2} + \frac{L}{6} \right) - 2 c_1 \log_2(N)\frac{\left(L-1\right)L}{2} \nonumber \\
		    & + c_2 \log_2^2(N) - c_2 \frac{\left(L-1\right)L}{2} + c_3 \log_2(N) + 1\nonumber \\
		    & = \frac{1}{2}c_1 \log_2^3(N) + \left(\frac{1}{2}c_1 + \frac{1}{2}c_2 \right)\log_2^2(N) + \left(\frac{1}{6}c_1 + \frac{1}{2}c_2 + c_3\right)\log_2(N) + 1. \label{appendix_eq7}
\end{align}

\bibliographystyle{siamplain}
\bibliography{references}

\begin{thebibliography}{10}

\bibitem{OpenBLAS}
{\em An optimized {BLAS} library}.
\newblock \url{http://www.openblas.net/}.
\newblock [Online; accessed December 12, 2018].

\bibitem{azad20163dmatmul}
{\sc A.~Azad, G.~Ballard, A.~Buluç, J.~Demmel, L.~Grigori, O.~Schwartz,
  S.~Toledo, and S.~Williams}, {\em Exploiting multiple levels of parallelism
  in sparse matrix-matrix multiplication}, SIAM J. Sci. Comput., 38 (2016),
  pp.~C624--C651, \url{https://doi.org/10.1137/15M104253X}.

\bibitem{communication_optimal_2}
{\sc G.~Ballard, A.~Bulu\c{c}, J.~Demmel, L.~Grigori, B.~Lipshitz, O.~Schwartz,
  and S.~Toledo}, {\em Communication optimal parallel multiplication of sparse
  random matrices}, in Proceedings of the Twenty-fifth Annual ACM Symposium on
  Parallelism in Algorithms and Architectures, SPAA '13, New York, NY, USA,
  2013, ACM, pp.~222--231.

\bibitem{benzi2000robust}
{\sc M.~Benzi, J.~K. Cullum, and M.~Tuma}, {\em Robust approximate inverse
  preconditioning for the conjugate gradient method}, SIAM J. Sci. Comput., 22
  (2000), pp.~1318--1332.

\bibitem{benzi2001stabilized}
{\sc M.~Benzi, R.~Kouhia, and M.~Tuma}, {\em Stabilized and block approximate
  inverse preconditioners for problems in solid and structural mechanics},
  Comput. Method. Appl. M., 190 (2001), pp.~6533--6554.

\bibitem{benzi1996sparse}
{\sc M.~Benzi, C.~D. Meyer, and M.~Tuma}, {\em A sparse approximate inverse
  preconditioner for the conjugate gradient method}, SIAM J. Sci. Comput., 17
  (1996), pp.~1135--1149.

\bibitem{Borstnik2014}
{\sc U.~Bor\v{s}tnik, J.~VandeVondele, V.~Weber, and J.~Hutter}, {\em Sparse
  matrix multiplication: The distributed block-compressed sparse row library},
  Parallel Comput., 40 (2014), pp.~47 -- 58.

\bibitem{Bowler_2012}
{\sc D.~R. Bowler and T.~Miyazaki}, {\em {O(N)} methods in electronic structure
  calculations}, Rep. Prog. Phys., 75 (2012), p.~036503.

\bibitem{BulucGilbert2012}
{\sc A.~Bulu\c{c} and J.~R. Gilbert}, {\em Parallel sparse matrix-matrix
  multiplication and indexing: Implementation and experiments}, SIAM J. Sci.
  Comput., 34 (2012), pp.~C170--C191.

\bibitem{Challacombe1999}
{\sc M.~Challacombe}, {\em A simplified density matrix minimization for linear
  scaling self-consistent field theory}, J. Chem. Phys., 110 (1999),
  pp.~2332--2342, \url{https://doi.org/10.1063/1.477969}.

\bibitem{dawson_ntpoly_2018}
{\sc W.~Dawson and T.~Nakajima}, {\em Massively parallel sparse matrix function
  calculations with {NTP}oly}, Comput. Phys. Commun., 225 (2018), pp.~154 --
  165.

\bibitem{duff1983multifrontal}
{\sc I.~S. Duff and J.~K. Reid}, {\em The multifrontal solution of indefinite
  sparse symmetric linear equations}, ACM T. Math. Software, 9 (1983),
  pp.~302--325.

\bibitem{gershgorin1931uber}
{\sc S.~A. Gershgorin}, {\em {\"U}ber die {A}bgrenzung der {E}igenwerte einer
  {M}atrix}, Izv. Akad. Nauk S.S.S.R., 6 (1931), pp.~749--754.

\bibitem{higham1997stable}
{\sc N.~J. Higham}, {\em Stable iterations for the matrix square root}, Numer.
  Algorithms, 15 (1997), pp.~227--242.

\bibitem{higham2008functions}
{\sc N.~J. Higham}, {\em Functions of matrices: theory and computation}, SIAM,
  2008.

\bibitem{jansik2007linear}
{\sc B.~Jans{\'\i}k, S.~H{\o}st, P.~J{\o}rgensen, J.~Olsen, and T.~Helgaker},
  {\em Linear-scaling symmetric square-root decomposition of the overlap
  matrix}, J. Chem. Phys., 126 (2007), p.~124104.

\bibitem{kohn1965self}
{\sc W.~Kohn and L.~J. Sham}, {\em Self-consistent equations including exchange
  and correlation effects}, Phys. Rev., 140 (1965), p.~A1133.

\bibitem{kruchinina2016parameterless}
{\sc A.~Kruchinina, E.~Rudberg, and E.~H. Rubensson}, {\em Parameterless
  stopping criteria for recursive density matrix expansions}, J. Chem. Theory
  Comput., 12 (2016), pp.~5788--5802.

\bibitem{levine2009quantum}
{\sc I.~N. Levine, D.~H. Busch, and H.~Shull}, {\em Quantum chemistry}, vol.~6,
  Pearson Prentice Hall Upper Saddle River, NJ, 2009.

\bibitem{Millam1997}
{\sc J.~M. Millam and G.~E. Scuseria}, {\em Linear scaling conjugate gradient
  density matrix search as an alternative to diagonalization for first
  principles electronic structure calculations}, J. Chem. Phys., 106 (1997),
  pp.~5569--5577, \url{https://doi.org/10.1063/1.473579}.

\bibitem{niklasson2004iterative}
{\sc A.~M. Niklasson}, {\em Iterative refinement method for the approximate
  factorization of a matrix inverse}, Phys. Rev. B, 70 (2004), p.~193102.

\bibitem{localized_inverse_factorization}
{\sc E.~H. Rubensson, A.~G. Artemov, A.~Kruchinina, and E.~Rudberg}, {\em
  Localized inverse factorization},  (2018),
  \url{https://arxiv.org/abs/1812.04919}.

\bibitem{rubensson2008recursive}
{\sc E.~H. Rubensson, N.~Bock, E.~Holmstr{\"o}m, and A.~M. Niklasson}, {\em
  Recursive inverse factorization}, J. Chem. Phys., 128 (2008), p.~104105.

\bibitem{CHT-PARCO-2014}
{\sc E.~H. Rubensson and E.~Rudberg}, {\em Chunks and {T}asks: A programming
  model for parallelization of dynamic algorithms}, Parallel Comput., 40
  (2014), pp.~328--343, \url{https://doi.org/10.1016/j.parco.2013.09.006}.

\bibitem{rubensson2016locality}
{\sc E.~H. Rubensson and E.~Rudberg}, {\em Locality-aware parallel block-sparse
  matrix-matrix multiplication using the chunks and tasks programming model},
  Parallel Comput., 57 (2016), pp.~87--106.

\bibitem{rubensson2007hierarchic}
{\sc E.~H. Rubensson, E.~Rudberg, and P.~Sa{\l}ek}, {\em A hierarchic sparse
  matrix data structure for large-scale {H}artree-{F}ock/{K}ohn-{S}ham
  calculations}, J. Comput. Chem., 28 (2007), pp.~2531--2537.

\bibitem{Ergo-JCTC-2011}
{\sc E.~Rudberg, E.~H. Rubensson, and P.~Sa{\l}ek}, {\em {Kohn-Sham} density
  functional theory electronic structure calculations with linearly scaling
  computational time and memory usage}, J. Chem. Theory Comput., 7 (2011),
  pp.~340--350, \url{https://doi.org/10.1021/ct100611z}.

\bibitem{Ergo-SoftwareX-2018}
{\sc E.~Rudberg, E.~H. Rubensson, P.~Sa{\l}ek, and A.~Kruchinina}, {\em Ergo:
  An open-source program for linear-scaling electronic structure calculations},
  SoftwareX, 7 (2018), pp.~107--111.

\bibitem{VandeVondele_2012}
{\sc J.~VandeVondele, U.~Bor\v{s}tnik, and J.~Hutter}, {\em Linear scaling
  self-consistent field calculations with millions of atoms in the condensed
  phase}, J. Chem. Theory Comput., 8 (2012), pp.~3565--3573.

\bibitem{Xiang2007}
{\sc H.~J. Xiang, J.~Yang, J.~G. Hou, and Q.~Zhu}, {\em Linear scaling
  calculation of band edge states and doped semiconductors}, J. Chem. Phys.,
  126 (2007), p.~244707, \url{https://doi.org/10.1063/1.2746322}.

\end{thebibliography}

\end{document}